\documentclass[11pt,leqno]{article}

\usepackage{amsmath,amsfonts,amscd,amssymb,theorem}

\long\def\comment#1\endcomment{}

\comment
\endcomment

\makeatletter
\begingroup
\gdef\th@dotted{\normalfont\itshape
  \def\@begintheorem##1##2{%
        \item[\hskip\labelsep \theorem@headerfont ##1\ ##2.]}%
\def\@opargbegintheorem##1##2##3{%
   \item[\hskip\labelsep \theorem@headerfont ##1\ ##2\ (##3).]}}
\endgroup
\makeatother

\theoremstyle{dotted}

\newtheorem{theorem}{Theorem}[section]
\newtheorem{lemma}[theorem]{Lemma}

\newtheorem{prop}[theorem]{Proposition}
\newtheorem{corr}[theorem]{Corollary}

\makeatletter
\begingroup
\gdef\th@upshape{\normalfont
  \def\@begintheorem##1##2{%
        \item[\hskip\labelsep \theorem@headerfont ##1\ ##2.]}%
\def\@opargbegintheorem##1##2##3{%
   \item[\hskip\labelsep \theorem@headerfont ##1\ ##2\ (##3).]}}
\endgroup
\makeatother

\theoremstyle{upshape}

\newtheorem{defn}[theorem]{Definition}
\newtheorem{remark}[theorem]{Remark}
\newtheorem{exa}[theorem]{Example}

\makeatletter
\renewcommand{\subsection}{\@startsection{subsection}{2}{0pt}{-3ex
plus -1ex minus -0.2ex}{-2mm plus -0pt minus
-2pt}{\normalfont\bfseries}} 
\renewcommand{\subsubsection}{\@startsection{subsubsection}{3}{0pt}{-3ex
plus -1ex minus -0.2ex}{-2mm plus -0pt minus
-2pt}{\normalfont\bfseries}} 
\makeatother

\makeatletter
\@addtoreset{equation}{section}
\makeatother

\newcommand{\cntrct}                
{\hspace{2pt}\raisebox{1pt}{\text{$\lrcorner$}}\hspace{2pt}}

\newcommand{\proof}[1][Proof.]{\smallskip\noindent{\em #1}}
\def\endproof{\hfill\ensuremath{\square}\par\medskip}

\renewcommand{\labelenumi}{{\normalfont(\roman{enumi})}}

\def\eqref#1{\thetag{\ref{#1}}}

\let\latexref=\ref
\def\ref#1{{\normalfont{\latexref{#1}}}}

\newcommand{\wt}{\widetilde}

\newcommand{\dg}{\dagger}

\newcommand{\idot}{{\:\raisebox{1pt}{\text{\circle*{1.5}}}}}
\newcommand{\hdot}{{\:\raisebox{3pt}{\text{\circle*{1.5}}}}}
\newcommand{\Ind}{\operatorname{Ind}}
\newcommand{\Pro}{\operatorname{Pro}}

\newcommand{\ppt}{\operatorname{\sf pt}}

\newcommand{\Z}{\mathbb{Z}}
\newcommand{\C}{\mathcal{C}}

\newcommand{\E}{\mathcal{E}}
\newcommand{\B}{\mathcal{B}}
\newcommand{\T}{\mathcal{T}}

\newcommand{\Hhom}{\operatorname{\mathcal{H}{\it om}}}

\newcommand{\Ext}{\operatorname{Ext}}
\newcommand{\Id}{\operatorname{Id}}
\newcommand{\Hom}{\operatorname{Hom}}
\newcommand{\RHom}{\operatorname{RHom}}

\newcommand{\Fun}{\operatorname{Fun}}

\newcommand{\V}{\operatorname{\sf V}}

\newcommand{\D}{\mathcal{D}}

\newcommand{\eps}{\epsilon}

\newcommand{\A}{\mathcal{A}}

\newcommand{\Ab}{\operatorname{Ab}}

\newcommand{\colim}{\operatorname{\sf colim}}
\newcommand{\hocolim}{\operatorname{\sf hocolim}}
\newcommand{\holim}{\operatorname{\sf holim}}
\renewcommand{\lim}{\operatorname{\sf lim}}

\newcommand{\Sets}{\operatorname{Sets}}
\newcommand{\Shv}{\operatorname{Shv}}

\newcommand{\id}{\operatorname{\sf id}}
\newcommand{\ev}{\operatorname{\sf ev}}

\newcommand{\wc}{\check}

\newcommand{\Y}{\operatorname{\sf Y}}

\newcommand{\cosk}{\operatorname{\sf cosk}}

\newcommand{\Ker}{\operatorname{Ker}}
\newcommand{\Coker}{\operatorname{Coker}}

\newcommand{\Cov}{\operatorname{Cov}}
\newcommand{\HCov}{\operatorname{HCov}}

\newcommand{\Pos}{\operatorname{Pos}}

\def\phi{\varphi}
\def\eps{\varepsilon}

\renewcommand{\SS}{{\sf S}}
\newcommand{\amod}{{\text{\rm -mod}}}

\newcommand{\Epi}{\operatorname{\sf Epi}}

\newcommand{\Mor}{\operatorname{Mor}}
\newcommand{\DMor}{\operatorname{DMor}}

\newcommand{\copr}{\sqcup}

\newcommand{\Ho}{\operatorname{Ho}}

\newcommand{\Tot}{\operatorname{\sf Tot}}

\newcommand{\I}{\operatorname{\sf I}}
\newcommand{\Ll}{\operatorname{\sf L}}
\newcommand{\Rr}{\operatorname{\sf R}}

\newcommand{\Dh}{\operatorname{\mathcal{D}\mathcal{H}}}
\newcommand{\Dd}{\operatorname{\sf D}}

\newcommand{\Cc}{\operatorname{\sf C}}

\title{What do abelian categories form?}

\author{D. Kaledin\thanks{Supported by the Russian Science
    Foundation, grant 21-11-00153.}}

\begin{document}

\maketitle

\tableofcontents

\section*{Introduction.}

Long time ago, when mathematics was ``science about numbers'' --
which by the way it still is in the Eastern languages, at least
nominally -- objects of mathematical study were typically elements
in a set. These days, they usually form a category. However, what
happens when these objects are themselves categories? Just small
categories without any adornment of course form a $2$-category, and
this is it: we have the category of functors $\Fun(I,I')$ for any
small categories $I$, $I'$, these are equipped with composition
functors, identity objects, and associativity and unitality
isomorphisms, and this is the end of the story. But things are much
less clear when we are talking about categories with additional
structures and/or categories of a special type.

The case in point is modern ``non-commutative'' or ``categorical''
algebraic geometry as formulated e.g.\ by Konsevich and Soibelman
\cite{ks}. As a slogan, this is the ``geometry of derived
categories'': one studies an algebraic variety $X$ by looking at its
derived category $\D(X)$ of coherent sheaves. But $\D(X)$ is more
than just a category. At the very least, it carries a triangulated
structure in the sense of Verdier \cite{verd}, but it is well-known
that this is not enough --- in particular, triangulated functors
between triangulated categories do not form a triangulated category,
and worse then that, while it is more-or-less clear what should be
this ``triangulated category of functors'' from $\D(X)$ to $\D(X')$,
it is not possible to recover it just from $\D(X)$ and $\D(X')$. The
correct object of study is a triangulated category ``with an
enhancement'', and the meaning of ``enhancement'' is a matter of
choice.

One rather radical choice that is becoming more popular is to take a
step back and say that in fact all categories should be equipped
with an enhancement. Triangulated categories then correspond to
stable enhanced categories, and stability is a condition and not a
structure. Stable functors between stable categories do form a
stable category, so the theory looks reasonably complete and
natural. However, in practice, in all the existing formalisms such
as e.g.\ quasicategories or complete Segal spaces, an enhanced
category is something rather large and dependent on arbitrary
choices, and it only makes sense to consider it up to a ``weak
equivalence'' of some sort. Thus to work with enhanced categories,
one has to use the cumbersome machinery of abstract homotopy theory
-- model categories, simplicial homotopy theory, and so on. This is
fine in topology, where this machinery is needed in any case, but
feels excessive in more algebraic applications, where people are
used to the simple and powerful homological algebra of \cite{toho},
and strongly prefer chain complexes to simplicial sets.

Perhaps for this reason, the most common technology of enhancements
used in categorical geometry is based on the notion of a {\em
  differential-graded} or {\em DG-category}; we refer the reader to
\cite{kel} for a very good overview. In this context, the question
raised in the title of this paper was addressed in the pioneering
paper \cite{tam}. Among other things, Tamarkin constructs the
correct DG-category of DG-functors between two DG-categories over a
fixed field $k$, and then studies in detail the composition functors
and all sorts of higher structures that arise in the theory. In
particular, if we restrict our attention to a single DG-category
$A_\idot$, then its {\em Hochschild Cohomology} $HH^\hdot(A_\idot)$
is defined as the algebra $\RHom^\hdot(\Id,\Id)$, where $\Id$ is the
identity endofunctor of $A_\idot$, and it carries an additional
structure of an $E_2$-algebra, or equivalently, of a
$B_\infty$-algebra, see \cite[Section 5.4]{kel}. This is crucial for
developing deformation theory of $A_\idot$, with first-order
deformations described by classes in $HH^2(A_\idot)$.

\medskip

Our goal in this paper is to some extent complementary to what was
done in \cite{tam}. Namely, we observe that if not all then at least
some triangulated categories that appear in geometry not only come
from DG-categories, but are also derived categories of something
abelian. So, assume that we consider derived categories $\D(\C)$,
$\D(\C')$ of some abelian categories $\C$, $\C'$. Can we recover the
``correct'' category of functors from $\D(\C)$ to $\D(\C')$ if we
remember not only $\D(\C)$, $\D(\C')$, but also the abelian
categories $\C$, $\C'$? If so, do we have a convenient model for
this category of functors, in terms of $\C$ and $\C'$?

\medskip

In a sense, this looks like a toy model for the whole theory, since
derived categories are somewhat special, and remembering the abelian
category $\C$ is even more restrictive. One advantage of this toy
model is that the resulting theory is absolute: while one can
consider abelian categories linear over a fixed field $k$, one can
also work without fixing the ground field. In fact, if one does fix
$k$, then the answer to our question has been known for a long time:
if $\C$ is small and the target category is large enough --- say,
the ind-completion $\Ind(\C')$ of a small abelian category $\C'$ ---
then $k$-linear left-exact functors $\C \to \Ind(\C')$ form an
abelian category (this essentially goes back to the famous
Gabriel-Popescu Theorem from 1960-es, see Example~\ref{GP0.exa} below
for a precise statement). Its derived category is exactly what we
want. In particular, it gives the correct Hochschild Cohomology
$HH^\hdot(\C) = \RHom^\hdot(\Id,\Id)$ that can be identified with
the DG version, and deformation theory in this context was
successfully developped in \cite{LV} and subsequent work.

Alternatively, instead of functors $\C \to \Ind(\C')$ one can
consider functors $\Ind(\C) \to \Ind(\C')$ that are continuous ---
that is, commute with filtered colimits. This gives the same
category, but it can be defined in larger generality: instead of
ind-completions of small abelian categories, one can consider
arbitrary finitely presentable abelian categories (see
Subsection~\ref{ind.subs} for more details).

In the absolute case, things are less well-studied, or at least,
less easy to find in the literature. If we look at deformation
theory, then the basic example of a first-order deformation that is
not linear over a field is the square-zero extension $\Z/p^2\Z$ of
the prime field $\Z/p\Z$. In this case, to get a deformation class,
one has to replace Hochschild Cohomology with the so-called Mac Lane
Cohomology $HM^\hdot(\Z/p\Z)$. One attempt to extend this to general
abelian categories was \cite{KL}, where a bunch of functor
categories are constructed, together with associated versions of
Hochschild Cohomology, and there are comparison theorems that show
how to recover Mac Lane Cohomology and some other generalizations of
Hochschild Cohomology. However, the emphasis in \cite{KL} was on
these cohomology theories and comparison maps. What we want to do
here is to concentrate on the $2$-categorical structure. For any
finitely presentable abelian categories $\C$, $\C'$, we construct
one particular functor category $\Mor(\C,\C')$, with its derived
version $\DMor_{st}(\C,\C')$, and we then construct the composition
functors, and show that objects in $\Mor(\C,\C')$
resp.\ $\DMor_{st}(\C,\C')$ indeed act naturally on $\C$
resp.\ $\D(\C)$.

\medskip

The basic idea of the construction is somewhat surprising but very
old; it again essentially goes back to Gabriel-Popescu
Theorem. Given an additive functor $E:\C \to \C'$ between abelian
categories, one reinterprets the condition that $E$ is left-exact as
a condition of $E$ being a sheaf for an appropriate Grothendieck
topology on $\C^o$ (sometimes called ``single-epi topology''). Then
one observes that one can drop additivity: being a sheaf and being
additive are independent conditions, one makes perfect sense without
the other. We then take finitely presentable abelian categories
$\C$, $\C'$, and consider the category $\Fun_c(\C,\C')$ of all
continuous functors $\C \to \C'$. The category $\Fun_c(\C,\C')$ is
abelian, and our $\Mor(\C,\C') \subset \Fun_c(\C,\C')$ is the full
subcategory formed by sheaves. This category is also abelian, we
have the left-exact fully faithful embedding $e:\Mor(\C,\C') \to
\Fun_c(\C,\C')$, and its left-adjoint associated sheaf functor
$a:\Fun_c(\C,\C') \to \Mor(\C,\C')$ is exact.

To extend this to derived categories, we start with the positive
part $\D^{\geq 0}(\C)$ of the derived category $\D(\C)$, and we use
the classic extension technology due to Dold \cite{dold}. Namely,
for any abelian $\A$, the Dold-Kan equivalence identifies the
category $\Fun(\Delta,\A)$ of cosimplicial objects in $\A$ with the
category $C^{\geq 0}(A)$ of chain complexes in $\A$ concentrated in
non-negative cohomological degrees. A functor $E:\A \to \A'$ to some
abelian $\A'$ then extends to a functor $\Dd(E):C^{\geq 0}(\A) \to
C^{\geq 0}(\A')$ by passing to cosimplicial objects and applying $E$
pointwise. More generally, if we have a functor $E:\A \to C^{\geq
  0}(\A') \cong \Fun(\Delta,\A')$, we can define its Dold extension
$\Dd(E):C^{\geq 0}(\A) \to C^{\geq 0}(\A')$ by applying $E$
pointwise and then restricting the resulting bisimplicial object to
$\Delta \subset \Delta \times \Delta$. If $\A = \C$, $\A' = \C'$ are
finitely presentable, then this construction sends continuous
functors to continuous functors, and it always sends pointwise
quasiisomorphisms to pointwise quasiisomorphisms, thus descends to a
functor
$$
\Dd:\D^{\geq 0}(\Fun_c(\C,\C')) \to \D^{\geq 0}(\Fun_c(C^{\geq
  0}(\C),\C')),
$$
where we identify $C^{\geq 0}(\Fun(-,-)) \cong \Fun(-,C^{\geq
  0}(-))$ and localize with respect to pointwise quasiisomorphisms.

We then denote by $\DMor(\C,\C')$ the derived category of our
functor category $\Mor(\C,\C')$, with $\DMor^{\geq 0}(\C,\C')$
standing for its positive part, and we observe that the derived
functor $R^\hdot e$ of the embedding $e$ provides a full embedding
$R^\hdot e:\DMor^{\geq 0}(\C,\C') \to \D^{\geq
  0}(\Fun_c(\C,\C'))$. On the other hand, say that a functor
$C^{\geq 0}(\C) \to C^{\geq 0}(\C')$ is {\em homotopical} if it
sends quasiisomorphisms to quasiisomorphisms, and let $\Dh^{\geq
  0}(\C,\C') \subset \D^{\geq 0}(\Fun_c(C^{\geq 0}(\C),\C'))$ be the
full subcategory formed by homotopical continuous functors. With
this notation, we prove the following (for precise statements, see
Theorem~\ref{dold.thm} and Corollary~\ref{dmor.corr}):
\begin{itemize}
\item An object $E \in \D^{\geq 0}(\Fun_c(\C,\C'))$ lies
  in the essential image of the embedding $R^\hdot e$ iff its Dold
  extension $\Dd(E)$ is homotopical. Moreover, $\Dd \circ
  R^\hdot(e):\DMor^{\geq 0}(\C,\C') \to \Dh^{\geq 0}(\C,\C')$ is an
  equivalence.
\end{itemize}
It is interesting to note that Dold himself used his extension
procedure slightly differently. He discovered that for any $E:\A \to
\A'$ whatsoever, $\Dd(E)$ sends chain-homotopic maps between
complexes to chain-homotopic maps (in our language, this is
Lemma~\ref{ch.le}), and then defined the derived functors $D^i(E)$,
$i \geq 0$ by applying $\Dd(E)$ to an injective resolution of some
$A \in \A$, and taking homology of the resulting complex. This
actually fits together quite nicely with the sheaf-theoretic
approach --- for any $E$, $D^0(E)$ is a sheaf, and $E \cong D^0(E)$
if and only if $E$ was a sheaf. The higher derived functors
$D^\hdot(E)$ are the homology objects of $R^\hdot e(a(E))$ (this is
Proposition~\ref{dold.prop}).

To extend our functors to the full derived category $\D(\C)$, we
need to know that the Dold extension $\Dd(E)$ commutes with
homological shifts; however, in general, this is not true and should
not be true. The reason for this is the additivity condition that we
dropped, and we now reinstate it on the derived level. This can be
done in several equivalent ways, see Proposition~\ref{st.prop}, but
the simplest is just to require that $E:\C \to C^\hdot(\C')$ becomes
additive when we project to $\D(\C')$. This distinguishes a full
triangulated subcategory $\DMor_{st}(\C,\C') \subset \DMor(\C,\C')$
of {\em stable objects}, and these --- as soon as they are bounded
from below, for technical reasons --- act naturally by functors
$\D(\C) \to \D(\C')$. Note that $\DMor_{st}(\C,\C')$ inherits a
$t$-structure from $\DMor(\C,\C')$, and its heart
$\Mor_{st}(\C,\C')$ is simply the category of additive left-exact
continuous functors $\C \to \C'$, just as in the $k$-linear
situation. However, the whole category is not the derived category
of its heart. The difference appears already in degree $2$, and
includes Mac Lane Cohomology classes responsible for non-linear
deformations.

\medskip

To finish the introduction, let us give a section-by-section
overview of what we do throughout the paper, but before that, let us
mention things that we do {\em not} do:
\begin{itemize}
\item We do not prove that our $\DMor(\C,\C')$ is indeed the
  category of stable enhanced functors $\D\C) \to
  \D(\C')$ in the homotopically enhanced world; indeed, doing
  this would require us to pick a model for this enhanced world, and
  none are too appealing. However, given what we do prove, making
  this last step in any particular model should be a trivial
  exercise.
\item We consistently restrict ourselves to finitely presentable
  abelian categories -- such as inductive completions $\Ind(\C)$ of
  small abelian categories -- and do not explore weaker finiteness
  conditions. It should be possible to do something for more general
  abelian categories, but for illustration purposes, we stick to the
  simplest possible case.
\item We do not explore at all the higher structures on Hochschild
  Cohomology and its generalization, the main content of
  \cite{tam}. We believe that there is a very interesting story to
  explore here, and that looking at things at the level of abelian
  categories might clarify the general theory, but this should be
  the subject of further research.
\item We do not touch deformation theory. We will return to this
  elsewhere.
\end{itemize}
When one reads this list of omissions, one realizes that very
little, possibly nothing in what remains is new. Thus the paper
should be treated as an overview, with the main goal of presenting
and maybe reassembling known things in a slightly different way, and
highlighting the main ideas. The ideas themselves are definitely not
new either, and not due to us (in particular, the idea of using
sheaves is borrowed from the exposition of the Gabriel-Popescu
Theorem in \cite{BD}, and the importance of dropping additivity is
inspired by the seminal paper \cite{JP} and subsequent work of
Pirashvili and others). All the proofs are presented for the sake of
completeness and for the convenience of the reader, and alternative
much earlier proofs are probably available in the literature. Having
said that, we now give an overview of the overview.

\subsection*{Overview of the paper.}

Section~\ref{gen.sec} contains various preliminaries. It should not
be understood as a self-contained introduction to category theory
and homological algebra; our goal is to fix notation, explain
precisely non-standard terminology that we use, and emphasize useful
things that are not usually emphasized (such as the fact that being
additive is a condition on a category and not a
structure). Subsection~\ref{cat.subs} is general category theory;
non-standard terminology here is {\em left-closed subcategories} and
{\em left-pointed categories}. Subsection~\ref{ind.subs} is
concerned with presentability and ind and pro-completions, in the
spirit of \cite{KS}. Subsection~\ref{ab.subs} deals with abelian
categories, and Subsection~\ref{hom.subs} is devoted to the derived
ones. We assume known all the standard material found in any
textbook on homological algebra, but we do discuss less standard
stuff such as homotopy limits and colimits, and relationship between
short exact sequences and bicartesian squares
(Remark~\ref{ab.sq.rem}), and then between distinguished triangles
and homotopy bicartesian squares (Example~\ref{D.exa}). The latter
is of course inspired by the notion of stability in the homotopy
enhanced world, but it is useful even in the unenhanced setting.

Section~\ref{top.sec} is devoted to Grothendieck topologies. We give
a very brief overview in Subsection~\ref{top.subs}, with some
illuminating examples such as topologies on a finite partially
ordered sets. Among other things, we recall that for any small
category $I$ equipped with a topology, and any finitely presentable
abelian category $\E$, the category $\Shv(I,\E)$ of $\E$-valued
sheaves on $I$ is abelian, and the embedding $e:\Shv(I,\E) \to
\Fun(I^o,\E)$ is left-exact, with an exact left-adjoint associated
sheaf functor $a$. Then in Subsection~\ref{cov.subs}, we turn to a
particular class of topologies --- namely, those generated by
coverings formed by a single morphism. To axiomatize the situation,
we introduce the notion of a {\em covering class} $F$, and we prove
one general result on existence of certain special coverings
(Lemma~\ref{hull.le}).

In Section~\ref{hc.sec}, we turn to hypercoverings. These are
usually understood as augmented simplicial sets of a certain type,
but in fact, a large part of the theory exists in much larger
generality --- namely, with $\Delta$ replaced by a more general
small category $I$ --- and proofs actually become easier when
unencumbered by simplicial combinatorics. Thus we take the liberty
of spending some type on hypercoverings of different types, for
different $I$. Subsection~\ref{pos.subs} deals with finite partially
ordered sets, Subsection~\ref{ml.subs} extends this further to a
class of categories that includes both $\Delta$ and the category
$\Pos$ of finite partially ordered sets, and then in
Subsection~\ref{del.subs}, we return to the standard simplicial
story. This is the only part of the paper where some of our results
might be new.

Both Section~\ref{top.sec} and Section~\ref{hc.sec} are completely
categorical; homological algebra first appears in
Section~\ref{dold.sec}. We start by recalling basic facts about the
Dold-Kan equivalence (we skip the proofs). We then introduce our
main character, the single-epi topology on an abelian category $\A$
(in the terminology of Subsection~\ref{cov.subs}, it corresponds to
the covering class of epimorphisms). We show that by virtue of the
Dold-Kan equivalence, hypercoverings in the single-epi topology can
be identified with left resolutions in the sense of homological
algebra. We then prove a general result,
Proposition~\ref{dold.prop}, that for any small category $I$ with a
covering class $F$, and any finitely presentable abelian category
$\E$, computes the derived functor $R^\hdot e$ of the embedding
$e:\Shv(I,\E) \to \Fun(I^o,\E)$ in terms of hypercoverings in
$I$. The result is in fact standard, and it also holds for more
general topologies --- ``local cohomology can be computed by
hypercoverings'' --- but the proof in our case is easy, and the end
result is nice: if $I=\A$, with the single-epi topology, then since
hypercoverings are left resolutions, the homology objects of
$R^\hdot e(a(E))$ for some functor $E:\A^o \to \E$ are simply the
Dold-derived functors of $E$.

Having finished with all the preliminaries, we can turn to the main
subject of the paper, namely, functor categories. In
Section~\ref{main.sec}, we define the category of functors
$\Mor(\C,\C')$ for any finitely presentable abelian categories $\C$,
$\C'$ and its derived category $\DMor(\C,\C')$, and we prove our
main extension results, Theorem~\ref{dold.thm} and
Corollary~\ref{dmor.corr}. These provide an action of the positive
part $\DMor^{\geq 0}(\C,\C')$ by functors $\D^{\geq 0}(\C) \to
\D^{\geq 0}(\C')$. Then in Section~\ref{stab.sec}, we introduce the
stability condition on functors, in several equivalent forms given
in Proposition~\ref{st.prop}, and show that the full subcategory
$\DMor_{st}^+(\C,\C') \subset \DMor(\C,\C')$ spanned by stable
objects bounded from below acts by functors $\D^+(\C) \to
\D^+(\C')$.

\subsection*{Acknowledgement.}

This paper grew out of a joint project with M. Booth and W. Lowen on
deformations of abelian categories, and I am grateful to both
colleagues for many fruitful discussions. I am also extremely
grateful to A. Efimov, V. Vologodsky and the anonymous referee for
quite a lot of useful remarks and comments, and a number of crucial
mistakes found in the first draft.

\section{Generalities.}\label{gen.sec}

\subsection{Categories and functors.}\label{cat.subs}

We denote by $\ppt$ the point category (one object, one morphism).
For any category $I$, we denote by $I^o$ the opposite category. For
any functor $\gamma:I \to \E$, we denote by $\gamma^o:I^o \to \E^o$
the opposite functor. For any object $E \in \E$, we denote by $E_I:I
\to \E$ the constant functor with value $E$. Somewhat
non-standardly, we will say that a full subcategory $I' \subset I$
is {\em left-closed} if for any map $f:i' \to i$ in $I$ with $i \in
I'$ we also have $i' \in I'$, and $I' \subset I$ is {\em
  right-closed} if ${I'}^o \subset I^o$ is left-closed. A category
$I$ is {\em connected} if any two objects are connected by a chain
of morphisms. For any object $e \in \E$, we denote by $I /_\gamma e$
the {\em left comma-fiber} of the functor $\gamma$, that is, the
category of pairs $\langle i,\alpha \rangle$, $i \in I$,
$\alpha:\gamma(i) \to e$ a morphism, we let $\sigma(e):I /_\gamma e
\to I$ be the forgetful functor, and we drop $\gamma$ from notation
when it is clear from the context (in particular, for any $i \in I$,
we shorten $I/_{\id} i$ to $I/i$). Dually, the {\em right
  comma-fiber} is $e \setminus_\gamma I = (I^o / e)^o$, and the {\em
  fiber} $I_e$ is the full subcategory $I_e \subset I/e$ spanned by
$\langle i,\alpha \rangle$ such that $\alpha$ is an isomorphism.

We treat a partially ordered set $J$ as a small category in the
usual way (objects are elements $j \in J$, there is a single
morphism from $j$ to $j'$ iff $j' \geq j$). A small category $I$ is
equivalent to a partially ordered set iff there is at most one
morphism between any two objects, and by abuse of terminology, we
will simply say that $I$ is a partially ordered set. Note that in
this case, for any $i \in I$, the comma-fiber $I / i$ is a
left-closed full subcategory $I / i \subset I$, and being full, it
is also a partially ordered set.

As another non-standard bit of terminology that will prove useful,
we say that a category $I$ is {\em left-pointed} if it has an
initial object $o \in I$, and $\{o\} \subset I$ is left-closed. Any
category $I$ can be turned into a left-pointed one by formally
adding an initial object $o$, and we will denote the resoluting
category by $I^<$. Any left-pointed category $I$ is of this type
(namely, we have $I \cong (I \setminus \{o\})^<$). The product $I_0
\times I_1$ of left-pointed categories $I_0$, $I_1$ is left-pointed,
and for any two categories $I_0$, $I_1$, we define their {\em
  extended product} by
\begin{equation}\label{st.eq}
I_0 * I_1 = (I^<_0 \times I^<_1) \setminus \{o \times o\},
\end{equation}
so that $I^<_0 \times I^<_1 \cong (I_0 * I_1)^<$.

\begin{exa}\label{V.exa}
The category $\ppt * \ppt$ can be naturally identified with the
opposite $\V^o$ to the partially ordered set $\V = \{0,1\}^<$ with
three elements $o$, $0$, $1$ and order relations $0,1 \geq o$.
\end{exa}

A functor between left-pointed categories is left-pointed if it
sends $o$ to $o$, and for any categories $I_0$, $I_1$, a functor
$\phi:I_0 \to I_1^<$ canonically extends to a left-pointed functor
$\phi^<:I_0^< \to I_1^<$. An {\em augmented functor} from a category
$I$ to some category $\E$ is a functor $E^<:I^< \to \E$, it is {\em
  $e$-aug\-mented} for some object $e \in \E$ if $E^<(o) = e$, and
an augmentation of a given functor $E:I \to \E$ is an augmented
functor $E^<:I^< \to \E$ equipped with an isomorphism $E^<|_I \cong
E$. Augmentations of a given functor $E$ form a category; by
definition, it has a terminal object $E^<$ iff the limit $\lim_IE$
exists, and we have $E^<(o) \cong \lim_IE$. In this case, we say
that the augmentation is {\em universal}. Giving an $e$-augmentation
of a functor $E:I \to \E$ is equivalent to giving a map $e_I \to E$
from the constant functor with value $e$. Dually, we denote $I^> =
(I^{o<})^o$, and a coaugmented functor $I \to \E$ is a functor $I^>
\to \E$; the universal coagumentation is given by the colimit
$\colim_I$ (if it exists). A category $\E$ is {\em complete}
resp.\ {\em finitely complete} if $\lim_IE$ exists for any small
resp.\ finite $I$ and functor $E:I \to \E$, and {\em cocomplete}
resp.\ {\em finitely cocomplete} if $\E^o$ is complete
resp.\ finitely complete. For example, the category $\Sets$ of all
sets is complete and cocomplete.

A {\em retract} of an object $i \in I$ in a category $I$ is an
object $i' \in I$ equipped with maps $a:i' \to i$, $b:i \to i'$ such
that $b \circ a = \id$. The composition $p=a \circ b:i \to i$ is
then idempotent, $p^2=p$, and $i'$ is the {\em image} of the
idempotent endomorphism $p$. The image is unique if it exists, and
is automatically preserved by any functor. A category is {\em
  Karoubi-closed} if every idempotent endomorphism of any object
admits an image. If a category $I$ is Karoubi-closed, then so is the
opposite category $I^o$. A category that is complete or cocomplete
is Karoubi-closed.

For any small category $I$ and arbitrary category $\E$, we let
$\Fun(I,\E)$ be the category of functors $E:I \to \E$. For any
functor $\gamma:I' \to I$ from a small $I'$, we let $\gamma^*E = E
\circ \gamma \in \Fun(I',\E)$. We have $\Fun(\ppt,\E)=\E$. For any
$E_0,E_1:I \to \E$, we denote by $\Hom_I(E_0,E_1)$ the set of maps
from $E_0$ to $E_1$, and we drop the index $I$ when it is clear from
the context. We also define a functor $\Hhom_I(E_0,E_1):I \to \Sets$
by
\begin{equation}\label{hhom.eq}
\Hhom_I(E_0,E_1)(i) = \Hom_{i \setminus
  I}(\sigma(i)^{o*}E_0,\sigma(i)^{o*}E_1),
\end{equation}
and we note that we have $\Hom_I(E_0,E_1) =
\lim_I\Hhom_I(E_0,E_1)$. For any $E \in \Fun(I',\E)$, the {\em
  left Kan extension} $\gamma_!E$ is a functor $\gamma_!E:I \to \E$
equipped with a map $E \to \gamma^*\gamma_!\E$ satisfying the usual
universal property. If $\colim_{I' / i}\sigma(i)^*E$ exists for any
$i \in I$, then the left Kan extension $\gamma_!E$ also exists, and
it is given by
\begin{equation}\label{kan.eq}
\gamma_!E(i) = \colim_{I' / i}\sigma(i)^*E, \qquad i \in I.
\end{equation}
If this happens for any $E$ --- for instance, if the target category
$\E$ is cocomplete --- then $\gamma_!:\Fun(I',\E) \to \Fun(I,\E)$ is
left-adjoint to $\gamma^*$. If $\gamma=\tau:I \to \ppt$ is the
tautological projection to a point, then $\tau_! = \colim_I$ is the
colimit itself. Dually, the {\em right Kan extension} $f_*E$ is
$(\gamma_!E^o)^o$, and there is a dual version of \eqref{kan.eq}
expressing $\gamma_*E$ in terms of limits over right comma-fibers.

\begin{remark}\label{univ.kan.rem}
The left Kan extension $\gamma_!E$ may exist even if some of the
colimits in \eqref{kan.eq} do not. Namely, for any functor $E':I \to
\E$ and object $i \in I$, a map $E \to \gamma^*E'$ induces a
$E'(i)$-coaugmentation of the functor $\sigma(i)^*E:I'/i \to \E$,
and let us say that $E'$ is a {\em universal left Kan extension} if
all these coaugmentations are universal. Then a universal left Kan
extension exists iff so do all the colimits in \eqref{kan.eq}, and a
universal left Kan extension is in particular a left Kan
extension. Under some assumptions -- for instance, if $I=\ppt$, or
if $\E$ has arbitrary products -- any left Kan extension is
universal, but in general, it is not true. For example, if $\E$ is a
discrete category with more than one object, then $\Fun(I,\E) \cong
\E$ iff a small category $I$ is connected, but since $\E$ has no
initial object, $\colim_IE$ never exists for a functor $E:I \to \E$
from an empty category $I$. Therefore for any functor $\gamma:I' \to
I$ between connected small categories, and any $E:I' \to \E$,
$\gamma_!E$ exists taulogically, but if at least one comma-fiber
$I'/i$ is empty, it is not universal.
\end{remark}

\begin{exa}\label{X.E.exa}
For any small $I$ and functor $X:I^o \to \Sets$, define the {\em
  category of elements} $IX$ as the category of pairs $\langle i,x
\rangle$, $i \in I$, $x \in X(i)$, with morphisms $\langle i,x
\rangle \to \langle i',x' \rangle$ given by morphisms $f:i \to i'$
such that $f(x') = x$. We then have the forgetful functor $\pi:IX \to
I$, $\langle i,x \rangle \mapsto i$, and for any functor $E:I^o \to
\E$ to a complete target category $\E$, we can define
\begin{equation}\label{hom.E}
\Hom(X,E) = \lim_{IX^o}\pi^{o*}E.
\end{equation}
Moreover, we can define a functor $\Hhom(X,E):I^o \to \E$ by
\begin{equation}\label{hhom.E}
\Hhom(X,E) = \pi^o_*\pi^{o*}E,
\end{equation}
and we then have $\Hom(E',\Hhom(X,E)) \cong \Hom(X,\Hhom(E',E))$ for
any $E' \in \Fun(I^o,\E)$, where $\Hhom(E',E)$ is given by
\eqref{hhom.eq}. If $\E = \Sets$, then \eqref{hom.E} simply computes
the set of morphisms $X \to E$, and \eqref{hhom.E} reduces to
\eqref{hhom.eq} by the dual version of \eqref{kan.eq}.
\end{exa}

For any integer $n \geq 0$, we denote by $[n]$ the ordinal
$\{0,\dots,n\}$ with the usual order, and when needed, we treat it
as a partially ordered set or a category. We have $[n]^o \cong [n]$
and $[n]^< \cong [n]^> \cong [n+1]$. For small $n$, $[0]=\ppt$ is a
point, and $[1]$ is the ``single arrow category'' that has two
objects $0$, $1$ and a single non-trivial arrow $0 \to 1$; functors
$[1] \to \E$ to some $\E$ correspond to arrows in $\E$. Functors
$[2] \to \E$ correspond to composable pairs of arrows $f$, $f'$; we
have the embeddings $s,t:[1] \to [2]$ onto the initial
resp.\ terminal segment of the ordinal $[2]$, and a functor $[2] \to
\E$ produces $f$ resp.\ $f'$ by restriction via $s$ resp.\ $t$. We
also have the embedding $m:[1] \to [2]$ onto $\{0,2\} \subset [2]$,
and restricting via $m$ produces the composition $f' \circ f$. More
generally, for any category $I$ and two functors $E_0,E_1:I \to \E$,
giving a map $f:E_0 \to E_1$ is equivalent to giving a functor
$\iota(f):[1] \times I \to \E$ whose restriction to $\{l\} \times
I$, $l=0,1$ is identified with $E_l$. A composable pair of maps $f$,
$f'$ gives a functor $\iota(f,f'):[2] \times I \to \E$ equipped with
isomorphisms $(s \times \id)^*\iota(f,f') \cong \iota(f)$, $(t
\times \id)^*\iota(f,f') \cong \iota(f')$, and we have a canonical
identification $\iota(f' \circ f) \cong (m \times
\id)^*\iota(f,f')$.

Analogously, commutative squares in a category $\E$ correspond to
functors $[1]^2 \to \E$, where $[1]^2 = [1] \times [1]$ is the
cartesian square of the single-arrow category $[1]$. The
single-arrow category $[1] \cong \ppt^<$ is left-pointed, and $\ppt
* \ppt \cong [1]^2 \setminus \{0 \times 0\}$ is the partially
ordered set $\V^o$ by Example~\ref{V.exa}, so that a commutative
square $\gamma:[1]^2 \to \E$ defines a augmented functor from $\V^o$
to $\E$. The square is cartesian iff the augmentation is universal
(that is, $\lim_{\V^o}\gamma$ exists, and the map $\gamma(0 \times
0) \to \lim_{\V^o}\gamma$ is an isomorphism). Dually, we have $[1]^o
\cong [1]$, so that we also have $[1]^2 \cong \V^>$, and a
commutative square is a coaugmented functor from $\V$; the square is
cocartesian iff the coaugmentation is universal.

\subsection{Inductive completions.}\label{ind.subs}

Let us recall some more advanced results on limits and inductive
completions that we will need (a good recent general reference for
all this is \cite[Chapter 6]{KS}). For a connected small category
$I$, with the tautological projection $\tau:I \to \ppt$, $\tau^*:\E
\to \Fun(I,\E)$ is fully faithful, so that by adjunction,
$\colim_IE_I \cong \lim_IE_I \cong E$ for any $E \in \E$ (and both
exist). A functor $\gamma:I' \to I$ is {\em cofinal} if $i \setminus
I'$ is connected for any $i \in I$. In this case, the dual version
of \eqref{kan.eq} shows that $\gamma_*E_{I'} \cong E_I$, $E \in \E$,
and that for any $E:I \to \E$, we have the adjunction isomorphism
\begin{equation}\label{cof.eq}
\colim_{I'}\gamma^*E \cong \colim_IE,
\end{equation}
where both sides exist at the same time. A category that has an
initial object is trivially connected, so that any functor that
admits a left-adjoint is cofinal. A useful example of such a
situation occures when $\gamma^o$ is a ``Grothendieck fibration'' of
\cite{SGA1} (see e.g.\ \cite[Section 1.3]{ka.adj} for a recent
overview with the same notation as here). In this case, the
embedding $I'_i \to I'/i$ admits a left-adjoint for any $i \in I$,
so that one can replace the left comma-fibers $I' / i$ in
\eqref{kan.eq} with the usual fibers $I'_i$.

\begin{defn}\label{filt.def}
A non-empty category $I$ is {\em directed} if \thetag{i} for any two
objects $i,i' \in I$, there exists an object $i'' \in I$ and maps $i
\to i''$, $i' \to i''$, and {\em filtered} if, moreover, \thetag{ii}
for any two maps $f,f':i \to i'$, there exists a map $g:i' \to i''$
such that $g \circ f = g \circ f'$.
\end{defn}

A directed category is obviously connected, so that if every right
comma-fiber $i \setminus I'$ of a functor $\gamma:I' \to I$ is
non-empty and directed for any $i \in I$, the functor $\gamma$ is
cofinal. If $I'$ is filtered, then all these comma-fibers satisfy
Definition~\ref{filt.def}~\thetag{ii} automatically, so they are
also filtered. If $I$ is a partially ordered set, then again,
Definition~\ref{filt.def}~\thetag{ii} is automatic, and $I$ is
filtered iff it is directed. A finitely cocomplete category $I$ is
trivially filtered. For any filtered $I$, the colimit functor
$\colim_I:\Fun(I,\Sets) \to \Sets$ preserves finite limits, and this
property is the main reason why the notion of a filtered category is
useful.

\begin{remark}
It is a pleasant exercise to check that the converse is also true --
if $\colim_I:\Fun(I,\Sets) \to \Sets$ preserves finite limits, then
$I$ is filtered.
\end{remark}

For any category $\C$, objects of the {\em inductive
  completion} $\Ind(\C)$ are pairs $\langle I,c \rangle$ of a small
filtered category $I$ and a functor $c:I \to \C$, and morphisms are
$$
\Hom(\langle I,c \rangle,\langle I',c' \rangle) = \lim_{i \in
  I^o}\colim_{i' \in I'}\Hom(c(i),c'(i')).
$$
We have the tautological full embedding $\iota:\C \to \Ind(\C)$, $c
\mapsto \langle \ppt,c \rangle$. The category $\Ind(\C)$ has
filtered colimits and is universal with this property. Namely, say
that a functor is {\em continuous} if it preserves filtered
colimits; then for any target category $\E$ that has filtered
colimits, and any functor $E:\C \to \E$, the left Kan extension
$\iota_!E:\Ind(\C) \to \E$ exists and is continuous, and is the
unique, up to a unique isomorphism continuous extension of $E$ to
$\Ind(\C)$. For any $C \in \Ind(\C)$ represented by a pair $\langle
I,c \rangle$, the natural projection $I \to \C / C$ is cofinal, so
that $\iota_!(E)(C) \cong \colim_{i \in I}E(c(i))$ by
\eqref{kan.eq} and \eqref{cof.eq}. Dually, the {\em projective
  completion} $\Pro(\C)$ is given by $\Pro(\C) =
\Ind(\C^o)^o$. Objects in $\Ind(\C)$ resp.\ $\Pro(\C)$ are also
called ind-objects resp.\ pro-objects in $\C$. If $\C = I$ is a
small category, then we have the Yoneda full embedding
\begin{equation}\label{yo.eq}
\Y:I \to \Fun(I^o,\Sets), \qquad \Y(i)(i') = \Hom(i',i),
\end{equation}
and $\iota_!\Y:\Ind(I) \to \Fun(I^o,\Sets)$ is also a full embedding
that identifies $\Ind(I)$ with the full subcategory in
$\Fun(I^o,\Sets)$ spanned by filtered colimits of representable
functors. Note that the Yoneda embedding \eqref{yo.eq}, hence also
the embedding $I \to \Ind(I)$ reflects monomorphisms (that is, a map
$f$ is a monomorphism in $I$ iff it is a monomorphism in $\Ind(I)$).

An object $c \in \C$ in a category $\C$ is {\em finitely
  presentable} or {\em compact} if the corepresentable functor
$\Hom(c,-)$ preserves filtered colimits. Let $\C_c \subset \C$ be
the full subcategory spanned by compact objects in a cocomplete
category $\C$; then for any full subcategory $I \subset \C_c$, the
full embedding $I \to \C_c \to \C$ canonically extends to a fully
faithful functor
\begin{equation}\label{ind.I.C}
\Ind(I) \to \Ind(\C_c) \to \C.
\end{equation}
A category $\C$ is {\em finitely presentable} if it is cocomplete,
and there exists a small full subcategory $I \subset \C_c$ such that
the functor \eqref{ind.I.C} is essentially surjective. Since it is
fully faithful, it is then automatically an equivalence. Moreover,
the functor $\Y_I:\C \to \Fun(I^o,\Sets)$ induced by the Yoneda
embedding \eqref{yo.eq} is a full embedding that preserves limits
and filtered colimits, so that in particular, filtered colimits in
$\C$ commute with finite limits, just as in the case $\C=\Sets$. For
any two finitely presentable categories $\C$, $\C'$, continuous
functors from $\C$ to $\C'$ form a well-defined category
$\Fun_c(\C,\C')$, and we have
\begin{equation}\label{fun.c.eq}
\Fun_c(\C,\C') \cong \Fun(I,\C'),
\end{equation}
where $I \subset \C_c \subset \C$ is a small subcategory such that
$\C \cong \Ind(I)$.

\begin{exa}\label{ind.exa}
Any colimit in a category $\C$ can be represented as a filtered
colimit of finite colimits. Therefore for any small finitely
cocomplete category $I$, the ind-completion $\Ind(I)$ is finitely
presentable. The Yoneda embedding $\Y_I$ induced by \eqref{yo.eq}
identifies $\Ind(I)$ with the full subcategory $\Fun_{ex}(I^o,\Sets)
\subset \Fun(I^o,\Sets)$ of functors $X:I^o \to \Sets$ that preserve
finite limits. Indeed, $\Fun_{ex}(I^o,\Sets)$ contains the
representable functors $\Y(i)$ and is closed under filtered
colimits, so that $\Ind(I) \subset \Fun_{ex}(I^o,\Sets)$, and on the
other hand, for any $X \in \Fun_{ex}(I^o,\Sets)$, the category of
elements $IX$ of Example~\ref{X.E.exa} is finitely cocomplete, thus
filtered, so that $X \cong \colim_{\langle i,x \rangle \in IX}Y(i)$
is in $\Ind(I)$.
\end{exa}

\begin{exa}\label{comp.exa}
The situation of Example~\ref{ind.exa} is in fact general. Namely,
for any small category $I$, an object $c \in \Ind(I)$ is compact iff
it is a retract of an object $i \in I \subset \Ind(I)$ (by
definition, $c \cong \colim_Ji$ for some filtered $J$ and functor
$i:J \to I$, and if $c$ is compact, the isomorphism $c \to
\colim_Ji$ must factor through $i(j) \in I$ for some $j \in
J$). Therefore for any finitely presentable category $\C \cong
\Ind(I)$, $I$ small, $\C_c \subset \C$ is essentially small, and we
then also have $\C \cong \Ind(\C_c) \cong \Fun_{ex}(\C_c^o,\Sets)$
(in particular, $\C$ is automatically complete). Since filtered
colimits of sets commute with finite limits, $\C_c \subset \C$ is
closed under finite colimits, thus finitely cocomplete.
\end{exa}

\subsection{Abelian categories.}\label{ab.subs}

A category $\C$ is {\em pointed} if it has an initial object $0$ and
a terminal object $1$, and the unique map $0 \to 1$ is an
isomorphism (so that $0$ is both an initial and a terminal object,
unique up to a unique isomorphism). For any two objects $A,B \in \C$
in a pointed category $\C$, we have a unique map $0:A \to B$ that
factors through $0$, so that the $\Hom$-sets $\Hom(-,-)$ are
naturally pointed. If a pointed category $\C$ admits finite products
and coproducts, then we have a natural map
\begin{equation}\label{copr.prod}
\begin{CD}
A \copr B @>{(\id \times 0) \copr (0 \times \id)}>> A \times B
\end{CD}
\end{equation}
for any $A,B \in \C$. Say that a category is {\em preadditive} if it
is pointed, has finite products and coproducts, and all the maps
\eqref{copr.prod} are isomorphisms. Thus we have $A \copr B \cong A
\times B$ canonically, and one denotes this object by $A \oplus B$
and calls it the {\em sum} of $A$ and $B$ (by extension, all
coproducts that exist in $\C$ are then called ``sums''). For any
preadditive category, $\Hom$-sets carry a structure of a commutative
monoid, with $0$ as the unity element, and compositions are monoid
maps. A category is {\em additive} if it is preadditive, and the
monoids $\Hom(A,B)$, $A,B \in \C$ are abelian groups (that is,
admits inverses). The category $\Ab$ of all abelian groups is
additive, and any additive category is automatically enriched over
$\Ab$ (that is, compositions are compatible with the abelian group
structures on $\Hom$-sets). The original reference for the notion of
an additive category is \cite{toho}; however, it is useful to
remember that the $\Ab$-enrichment required in \cite[Section
  1.3]{toho} is actually automatic and unique, so that being
additive is a condition on a category and not a structure. We also
note that the property of being additive is self-dual --- the
opposite $\C^o$ to an additive category $\C$ is also additive.

\begin{remark}\label{preadd.rem}
The term ``preadditive'' used above is non-standard (and sometimes
it appears in the literature with a different meaning). It seems
that there is no standard term.
\end{remark}

The notion of an abelian category also goes back to \cite{toho}, and
it is also self-dual. Namely, for any map $f:A \to B$ in a pointed
category $\A$, the {\em kernel} is given by $\Ker f = A \times_B 0$,
and dually, the {\em cokernel} is $\Coker f = (\Ker f^o)^o$. Both
need not exist in general; an additive category is abelian if it has
all kernels and cokernels (``axiom $AB1$''), and for any morphism $f:A
\to B$, with kernel $k:\Ker f \to A$ and cokernel $c:B \to \Coker
f$, the natural map $\Coker k \to \Ker c$ is an isomorphism
($AB2$). An additive category satisfying $AB1$ is finitely complete
and cocomplete, and if it satisfies $AB2$, it is also
Karoubi-closed. An {\em abelian subcategory} $\A' \subset \A$ in an
abelian category $\A$ is a full subcategory closed under finite
sums, kernels and cokernels (so that in particular, $\A'$ is also
abelian). Alternatively, a {\em short exact sequence} in a pointed
category $\A$ is a sequence
\begin{equation}\label{ab1.dia}
\begin{CD}
A @>{i}>> B @>{p}>> C
\end{CD}
\end{equation}
such that $p \circ i = 0$, $\Ker p$ and $\Coker i$ exist, and the
maps $A \to \Ker p$, $\Coker i \to C$ are isomorphisms; then
another way to phrase $AB2$ is to say that every map $f:A \to
B$ admits a decomposition
\begin{equation}\label{ab2.dia}
\begin{CD}
C_0 @>>> A @>>> C_1 @>>> B @>>> C_2
\end{CD}
\end{equation}
such that $C_0 \to A \to C_1$ and $C_1 \to B \to C_2$ are short
exact sequences (note that such a decomposition is necessarily
unique).

\begin{remark}\label{ab.sq.rem}
The following repackaging of the notion of a short exact sequence is
sometimes useful. Giving a sequence \eqref{ab1.dia} with $p \circ i
= 0$ is equivalent to giving a commutative square
\begin{equation}\label{ab1.sq}
\begin{CD}
A @>{i}>> B\\
@VVV @VV{p}V\\
0 @>>> C.
\end{CD}
\end{equation}
Then a sequence \eqref{ab1.dia} is exact on the left resp.\ on the
right iff the corresponding square \eqref{ab1.sq} is cartesian
(equivalently, $\Ker p$ exists and $A \to \Ker p$ is an isomorphism)
resp.\ cocartesian (equivalently, $\Coker i$ exists, and $\Coker i
\to C$ is an isomorphism). The sequence is exact if the square is
    {\em bicartesian} (that is, cartesian and cocartesian at the
    same time).
\end{remark}

Grothendieck lists further conditions of increasing strength that
one can impose on an abelian category $\C$ --- it can have arbitrary
coproducts ($AB3$), coproducts of short exact sequences can be exact
($AB4$), and the same can hold for filtered colimits (or
equivalently, filtered colimits can commute with finite limits --
this is $AB5$). There is also a further property $AB6$ whose real
importance has begun to emerge only recently, so we skip it. The
additional properties are not self-dual; one says that an abelian
category $\C$ satisfies $ABN^*$, $N=3,4,5,6$ if $\C^o$ satisfies
$ABN$. The category $\Ab$ satisfies $AB5$ and $AB4^*$ (there is a
theorem that an abelian category satisfying $AB5$ and $AB5^*$ is
trivial).

A {\em generator} of an abelian category $\C$ is an object $U \in
\C$ such that $\Hom(U,-)$ is faithful (equivalently, for any
morphism $f:M \to M'$ in $\C$, $\Hom(U,f)=0$ implies $f=0$). A {\em
  Grothendieck abelian category} is an abelian category $\C$
satisfying $AB5$ that admits a generator. One of the main results of
\cite{toho} is that a Grothendieck abelian category has enough
injectives (that is, any $A \in \C$ admits a monomorphism $A \to I$
with injective $I$).

One usually defines {\em additive functors} for additive categories,
but it is useful to do it slightly more generally. Say that a
functor $E:\A \to \B$ between categories with finite products is
{\em additive} if it commutes with finite products (that is, for any
$A,B \in \A$, the natural map $E(A \times B) \to E(A) \times E(B)$
is an isomorphism). This is again a condition and not a structure;
however, if $\A$ and $\B$ are additive, then an additive functor $E$
is automatically enriched over $\Ab$. On the other hand, if $\A$ is
additive and $\B = \Sets$, then $E$ automatically and uniquely
factors through the forgetful functor $\Ab \to \Sets$ (more
precisely, the forgetful functor $\Fun(\A,\Ab) \to \Fun(\A,\Sets)$
induces an equivalence between the full subcategories spanned by
additive functors). In particular, since filtered colimits of sets
commute with finite products, all objects $E \in \Ind(\A) \subset
\Fun(\A^o,\Sets)$ in the inductive completion of a small additive
category $\A$ are additive, $\Ind(\A)$ is an additive category, and
the full embedding $\Ind(\A) \subset \Fun(\A^o,\Sets)$ factors
through a full embedding $\Ind(\A) \subset \Fun(\A^o,\Ab)$.

A functor between abelian categories is {\em left-exact} resp.\ {\em
  right-exact} if it commutes with finite limits resp.\ finite
colimits; in particular, a left or right-exact functor is
automatically additive. Alternatively, an additive functor is left
resp.\ right-exact iff it sends short exact sequences to sequences
exact on the left resp.\ on the right (the simplest way to see that
the two notions are equivalent is to use the description of short
exact sequences in terms of squares \eqref{ab1.sq}). A functor is
{\em exact} if it is both left and right-exact.

If a small category $\A$ is abelian, then $\Ind(\A)$ is finitely
presentable by Example~\ref{ind.exa}, and it is a Grothendieck
abelian category (this is well-known but rather non-trivial, see
below Example~\ref{GP.exa}). Conversely, any finitely presentable
abelian category $\E$ satisfies $AB5$ almost by definition --
filtered colimits commute with finite limits and colimits -- and for
any small full subcategory $I \subset \E_c \subset \E$ such that $\E \cong
\Ind(I)$, the sum of all objects $i \in I \subset \E$ is a generator,
so that $\E$ is a Grothendieck abelian category. The full
subcategory $\E_c \subset \E$ of compact objects is additive, and by
Example~\ref{comp.exa}, it is essentially small and has cokernels
(and all finite colimits).

\begin{defn}\label{preab.def}
A small additive category $\A$ with cokernels is {\em preabelian} if
$\Ind(\A)$ is abelian.
\end{defn}

Any small abelian category is preabelian, but the converse is not
the case, so the notion is not vacuous -- for example, the category
$R\amod$ of left modules over a ring $R$ is abelian, $M \in R\amod$
is compact iff it is the cokernel of a map $f:R^n \to R^m$ for some
integers $m,n \geq 0$, so that $R\amod$ is obviously finitely
presentable, but $(R\amod)_c$ is abelian only if the ring $R$ is
left-coherent. As in Remark~\ref{preadd.rem}, the term
``preabelian'' is not standard, but it seems that there is no
standard term. An abstract characterization of preabelian categories
can be found in \cite{sch} (under the name ``ind-abelian''). By
Example~\ref{comp.exa}, a preabelian category $\A$ is of the form
$\E_c$, $\E$ abelian finitely presentable if and only if $\A$ is
Karoubi-closed.

\subsection{Derived categories.}\label{hom.subs}

We denote by $C^{\hdot}(\E)$ the category of chain complexes
$M^\hdot = \langle M^\hdot,d \rangle$ in an additive category $\E$,
and we let $C_\idot(\E)$ be the same category but with complexes
indexed by homological rather than cohomological degrees, with the
convenion being $M^i = M_{-i}$. The {\em homological shift}
$M^\hdot[n]$ of a complex $M^\hdot$ by an integer $n$ is given by
$(M^\hdot[n])^i = M^{i+n}$. The {\em cone} of a map $f:M^\hdot \to
N^\hdot$ in $C^\hdot(\E)$ is given by $\Cc(f)^i = N^i \oplus
M^{i+1}$, with the usual upper-triangular differential, and if $\E$
is abelian, we assume known the standard notions of homology
objects, an acyclic complex, a quasiisomorphism and so on. We denote
by $C^+(\E),C^{-}(\E) \subset C^\hdot(\E)$ the full subcategories
spanned by complexes bounded from below ($M^i = 0$ for $i \ll 0$)
resp.\ from above ($M^i = 0$ for $i \gg 0$), and we let
$C^\hdot_b(\E) = C^+(\E) \cap C^-(\E) \subset C^\hdot(\E)$ be the
full subcategory of bounded complexes. If $\E$ is finitely
presentable, with the subcategory $\E_c \subset \E$ of compact
objects, then $C^\hdot(\E)$ is also finitely presentable, and
$C^\hdot_b(\E_c) \subset \C^\hdot(\E)$ is the full subcategory of
compact objects. In particular, for a small abelian category $\A$,
we have $C^\hdot(\Ind(\A)) \cong \Ind(C^\hdot_b(\A))$.

The {\em localization} $h(\C,W)$ of a category $\C$ with respect to
a class of morphisms $W$ is a category $h(\C,W)$ equipped with a
functor $h:\C \to h(\C,W)$ that inverts all morphisms in $W$ and is
universal with this property: any functor $\C \to \E$ to some $\E$
that inverts all morphisms in $W$ factors through $h$, uniquely up
to a unique isomorphism. For any small category $I$, we denote by
$W^I$ the class of maps the functor category $\Fun(I,\C)$ that are
pointwise in $W$. If both localizations $h(\C,W)$,
$h(\Fun(I,\C),W^I)$ exist, we have the tautlogical functor $h(\C,W)
\to h(\Fun(I,\C),W^I)$, $c \mapsto c_I$, and the {\em homotopy
  limit} and {\em homotopy colimit} are by definition its left
resp.\ right-adjoint functors
\begin{equation}\label{holim.eq}
\hocolim_I,\holim_I:h(\Fun(I,\C),W^I) \to h(\C,W),
\end{equation}
if they exist. If $\holim_I$ exists, then an augmented functor
$c:I^< \to \C$ gives rise to a comparison map $c(o) \to \holim_Ic$,
and we say that the augmentation is {\em homotopy universal} if the
map is an isomorphism in $h(\C,W)$; in particular, a commutative
square $[1]^2 \to \C$ is {\em homotopy cartesian} if it is homotopy
universal when consider as an augmented functor (just as in the
non-homotopical case). Dually, if $\hocolim_I$ exists, we have the
notion of a homotopy universal coaugmentation and of a homotopy
cocartesian square.

Localization does not always exist (if $\C$ is large, there could be
set-theoretical issues), and is notoriously difficult to construct
explicitly and to describe (in particular, even if relevant
localizations exists, constructing homotopy limits and colimits is a
highly non-trivial task). One additional structure that helps to
control localization is that of a {\em model category} of Quillen,
see e.g.\ \cite{qui}, \cite{DS}, \cite{hov} (although \cite{hov} has
to be used carefully since the author takes the liberty of redefinng
standard notions according to his needs). Since we will only use it
tangentially (e.g.\ below in Lemma~\ref{hur.le}), we do not give any
details; let us just mention that if $\C$ is a model category and
$W$ is the class of weak equivalences, then by \cite{DHKS},
$h(\Fun(I,\C),W^I)$ exists for any finite $I$, and so do the
homotopy limit and colimit functors of \eqref{holim.eq}. In
particular, homotopy cartesian and homotopy cocartesian squares are
well-defined in any model category $\C$.

In the abelian context, the most common example of a localization is
the category $C^\hdot(\C)$ of chain complexes in an abelian category
$\C$ whose localization with respect to the class of
quasiisomorphisms produces the {\em derived category} $\D(\C)$. If
$\C$ has enough injectives, this localization can be constructed by
model category techniques, but there is a simpler and earlier
alternative that works under much milder restrictions: one first
constructs the {\em homotopy category} $\Ho(\C)$ of chain complexes
and chain-homotopy classes of maps between them, and then applies a
general localization theorem of Verdier \cite{verd} (this relies on
the structure of a {\em triangulated category} on $\Ho(\C)$ that was
introduced by Verdier specifically for this purpose). We refer the
reader to any standard textbook on homological algebra such as
\cite{W} or \cite{GM} for basic facts on derived categories; in
particular, we assume known the fact that the derived category
$\D(\C)$ is additive, and the dual notions of a total right-derived
functor $R^\hdot E:\D(\C) \to \D(\E)$ resp.\ left-derived functor
$L^\hdot E$ of a left-exact resp.\ right-exact functor $E:\C \to \E$
between abelian categories.

For any abelian category $\C$ and small category $I$, the functor
category $\Fun(I,\C)$ is abelian, and we simplify notation by
writing $\D(I,\C) = \D(\Fun(I,\C))$. If the abelian category $\C$ is
finitely presentable, then for any finitely presentable $I$, the
continuous functor category $\Fun_c(I,\C)$ is also abelian by
\eqref{fun.c.eq}, and we write $\D_c(I,\C) = \D(\Fun_c(I,\C)) \cong
\D(I_c,\C)$. Alternatively, $\D(I,\C)$ for a small $I$ and abelian
$\C$ can by obtained by localizing the functor category
$\Fun(I,C^\hdot(\C)) \cong C^\hdot(\Fun(I,\C))$ with respect to the
class of pointwise quasiisomorphisms, and if $I$ is finite, the
homotopy limit and colimit \eqref{holim.eq} both exists and are
given by $\hocolim_I = L^\hdot \colim_I$, $\holim_I = R^\hdot
\lim_I$. An object $E \in \D(I,\C)$ tautologically defines a functor
$\D(E):I \to \D(\C)$, so that we have a comparison functor
\begin{equation}\label{D.E}
\D:\D(I,\C) \to \Fun(I,\D(\C)).
\end{equation}
This functor is {\em not} an equivalence unless $I=\ppt$.

\begin{exa}\label{D.exa}
Take $I=[1]$, the single-arrow category. For any abelian $\C$,
objects in $\Fun([1],\C)$ are arrows in $\C$, and taking cokernel of
an arrow provides a right-exact functor $\Coker:\Fun([1],\C) \to \C$
with derived functor $L^\hdot\Coker:\D([1],\C) \to \D(\C)$. For any
$E \in \D([1],\C)$, $\D(E):[1] \to \D(\C)$ is an arrow in the
derived category $\D(\C)$, and $L^\hdot\Coker(E)$ gives its cone in
the sense of the triangulated structure on $\D(\C)$. However, this
version of the cone is functorial. The necessary rigidity is added
exactly by lifting $\D(E)$ to an object $E \in
\D([1],\C)$. Analogously, distinguished triangles in $\D(\C)$ can be
naturally rigidified by considering squares \eqref{ab1.sq} in
$\D([1]^2,\C)$ that are homotopy bicartesian. Every such square
produces a distinguished triangle after applying \eqref{D.E}, and
conversely, any distinguished triangle lifts to such a square; the
lifting is unique but only up to a non-unique isomorphism.
\end{exa}

One can also consider the category $C_{\geq 0}(\C) = C^{\leq 0}
\subset \C_\idot(\C) = C^\hdot(\C)$ of complexes concentrated in
non-negative homological (non-positive cohomological) degrees; its
localization produces the full subcategory $\D^{\leq 0}(\C) \subset
\D(\C)$ of {\em connective objects}, a part of a standard
$t$-structure on $\D(\C)$. Dually, localizing the category $C^{\geq
  0}(\C)$ of complexes concentrated in non-negative cohomological
degrees produced the full subcategory $\D^{\geq 0}(\C) \subset
\D(\C)$ of {\em coconnective objects}, another part of the standard
$t$-structure. For generalities on $t$-structures, see \cite{BBD};
let us just recall that the embedding $\D^{\leq 0}(\C) \subset
\D(\C)$ admits a right-adjoint canonical truncation functor
$\tau^{\leq 0}:\D(\C) \to \D^{\leq 0}(\C)$, and $\D^{\leq 0}(\C)
\cap \D^{\geq 0}(\C) \cong \C$, so that $\tau^{\geq 0}$ induces a
functor $\D^{\geq 0}(\C) \to \C$.

\section{Topologies and coverings.}\label{top.sec}

\subsection{Recollection on Grothendieck topologies.}\label{top.subs}

The original reference for Grothendieck topologies and topos theory
is \cite{SGA4}, and a very concise and useful overview can be found
in \cite[Chapter 0.3]{J}. Let us recall the basics of the story.

By definition, a {\em sieve} on an object $i \in I$ in a small
category $I$ is a subfunctor in the representable functor $\Y(i) =
\Hom(-,i):I^o \to \Sets$. The collection of all sieves on an object
$i$ is denoted $\Omega(i)$, and for any map $f:i' \to i$ and sieve
$s \in \Omega(i)$, $f^*s = s \times_{\Y(i)} \Y(i')$ is a sieve on
$i'$, so that $\Omega$ is itself a contravariant functor $I^o \to
\Sets$. A {\em Grothendieck topology} on $I$ is given by
collections of sieves $T(i)$, one for each object $i \in I$
satisfying the following axioms:
\begin{enumerate}
\item for any $i \in I$, the maximal sieve $\Y(i)$ is in $T(i)$,
\item for any map $f:i \to i'$ and $s \in T(i')$, $f^*s \in T(i)$
  (in other words, $T \subset \Omega$ is a subfunctor), and
\item if for some $s \in \Omega(i)$, $s' \in T(i)$ we have $f^*s \in
  T(i')$ for any $f:i' \to i$ in $s'(i') \subset \Hom(i',i)$,
  then $s \in T(i)$.
\end{enumerate}
For any $i \in I$, $\Omega(i)$ is a partially ordered set with
respect to the inclusion, and the axioms \thetag{i}-\thetag{iii}
imply that $T(i) \subset \Omega(i)$ is right-closed and closed under
intersections, so that $T(i)^o$ is a directed partially ordered set.

A functor $E:I^o \to \Sets$ is a {\em separable presheaf} resp.\ a
{\em sheaf} with respect to a topology $T$ if for any $i \in I$, $s
\in T(i)$, the map $E(i) \to \Hom(s,E)$ is injective
resp.\ bijective. If we compute $\Hom(s,E)$ by \eqref{hom.E}, then
this makes sense for a functor $E:I^o \to \E$ to any complete target
category $\E$, so that the notion of being a sheaf or a separable
presheaf is also defined for $\E$-valued functors. Explicitly, the
category of elements $Is$ of Example~\ref{X.E.exa} is equivalent to
the full subcategory $I /_s i \subset I / i$ spanned by arrows $f
\in s(i) \subset \Hom(i',i)$, and we have
\begin{equation}\label{hom.eq}
\Hom(s,E) = \lim_{i' \in (I  /_s i)^o}E(i').
\end{equation}
Denote by $\Shv(I,\E) \subset \Fun(I^o,\E)$ the full subcategory
spanned by sheaves. Then if the target category $\E$ is finitely
presentable --- thus in particular, complete and cocomplete --- the
embedding $\Shv(I,\E) \to \Fun(I^o,\E)$ admits a left-adjoint {\em
  associated sheaf functor} $a:\Fun(I^o,\E) \to \Shv(I,\E)$. To
construct it, define a functor $a_0:\Fun(I^o,\E) \to \Fun(I^o,\E)$
by
\begin{equation}\label{a.eq}
a_0(E)(i) = \colim_{s \in T(i)^o}\Hom(s,\E), \qquad i \in I, E \in
\Fun(I^o,\E).
\end{equation}
This is functorial in $i$ since so is $T(i)$ (slightly more
precisely, the partially ordered sets $T(i)$, $i \in I$ fit together
into a Grothendieck fibration $\gamma:\T \to I$ whose fibers are
partially ordered sets $T(i)$, with the inclusion order, and then
$a_0(E) = \gamma^o_!\Hom(-,E)$ is the left Kan extension along the
opposite functor $\gamma^o:\T^o \to I^o$). Then the natural maps
$E(i) \to \Hom(s,E)$ provide a functorial map $E \to a_0(E)$, and
one checks that $a_0(E)$ is a separated presheaf for any $E$, and a
sheaf if $E$ is separated (for $\E=\Sets$, this is \cite[Proposition
  3.2]{SGA4-verd}, and the general case reduces to this by the
Yoneda embedding \eqref{yo.eq}). Therefore $a_0^2(E) = a_0(a_0(E))$
is a sheaf for any $E$, and we have an isomorphism $a_0^2 \cong e
\circ a$ for a unique functor $a:\Fun(I^o,\E) \to \Shv(I,\E)$, while
the map $E \to a_0(E)$ provides the adjunction map $E \to
e(a(E))$. By adjunction, $a$ commutes with arbitrary colimits, and
since the colimits in \eqref{a.eq} are filtered, it also commutes
with finite limits. A functor $E:I^o \to \E$ is a sheaf iff the
adjunction map $f:E \to e(a(E))$ is an isomorphism, and in fact it
suffices to require that it admits a splitting $g:e(a(E)) \to E$, $g
\circ f = \id$ (because then $a(g) \circ a(f) = \id$, and since
$a(f)$ is invertible, $a(g \circ f) = a(g) \circ a(f) = \id:a(E) \to
a(E)$, so that $g \circ f = \id$ since $e$ is fully faithful).

\begin{exa}\label{shv.ab.exa}
If the target category $\E$ is abelian, then $\Shv(I,\E)$ is also
abelian, with kernels taken pointwise and cokernels created by the
associated sheaf functor $a$ (that is, for any map $f:E_0 \to E_1$
between $E_0,E_1 \in \Fun(I^o,\E)$ that are actually sheaves, $\Ker
f$ is a sheaf, and $a(\Coker f)$ is a cokernel in $\Shv(I,\E)$). To
check $AB2$, note that $a$ preserves cokernels by definition but also
commutes with finite limits, hence preserves kernels and short exact
sequences \eqref{ab1.dia}; then a decomposition \eqref{ab2.dia} in
$\Shv(I,\E)$ is obtained by applying $a$ to the corresponding
decomposition in $\Fun(I^o,\E)$. Moreover, $\Shv(I,\E)$ is a
Grothendieck abelian category that satisfies $AB3^*$ (where we recall
that $\E$ is finitely presentable by our standing assumption, to
insure the existence of the functor $a$). Indeed, products in
$\Shv(I,\E)$ are products in $\Fun(I^o,\E)$, and $a$ commutes with
colimits and finite limits, so that $AB5$ and $AB3^*$ are inherited
from $\Fun(I^o,\E)$; to obtain a generator, it suffices to take the
sum of all objects $a(\Y^i(E))$, $i \in I$, $E$ a fixed generator of
$\E$, $\Y^i(E)$ the corepresentable functor given by
\begin{equation}\label{yo.i.eq}
\Y^i(E)(i') = E[\Hom(i',i)], \qquad i' \in I,
\end{equation}
where the right-hand side is shorthand for ``sum of copies of
$E$ numbered by elements in the set $\Hom(i',i)$''. If $I$ has
finite coproducts, one can also consider the full subcategory
$\Shv_{add}(I,\E)$ spanned by sheaves that are additive; since
filtered colimits in $\E$ commute with finite products, the functor
$a$ preserves additivity, and $\Shv_{add}(I,\E) \subset \Shv(I,\E)$
is an abelian subcategory.
\end{exa}

\begin{remark}
In fact {\em any} Grothendieck abelian category satisfies $AB3^*$
but this is a rather non-trivial theorem; for the categories
$\Shv(I,\E)$ of Example~\ref{shv.ab.exa}, the claim is obvious.
\end{remark}

\begin{exa}\label{T.exa}
For any small $I$, the minimal topology $T_{min}$ consists of the
maximal sieves $\Y(i)$, $i \in I$; we have $T_{min} \cong \ppt_I$,
the functor sending everything to the one-point set $\ppt$, and the
corresponding embedding
\begin{equation}\label{min.t.eq}
1:\ppt_I \to \Omega
\end{equation}
is given by $1(i)(\ppt) = \Y(i) \in \Omega(i)$, $i \in I$ (note that
$\ppt_I$ is the terminal object in $\Fun(I^o,\Sets)$). Sheaves for
the minimal topology are all functors $I^o \to \E$. The maximal
topology $T_{max}$ consists of all sieves including the empty one,
$T_{max} \cong \Omega$, but it is not very interesting since the
only sheaf is the constant functor $I^o \to \E$ sending everything
to the terminal object in $\E$. However, if $I$ is left-pointed,
then all {\em non-empty} sieves also form a topology that we call
{\em submaximal}; for this topology, we have $\Shv(I^o,\E) \cong
\E$, with the equivalence given by evaluation at the initial object
$o \in I$.
\end{exa}

\begin{exa}\label{pos.top.exa}
Let $J$ be a partially ordered set. Then a sieve on an object $j \in
J$ is the same thing as a left-closed subset in the comma-set $J / j
\subset J$. For any subset $J' \subset J$ and element $j \in J$, let
$T_{J'}(j) \subset \Omega(j)$ consist of sieves $J_0 \subset J/j$
that contain $J' \cap (J /j)$. Then $T_{J'}$ is a Grothendieck
topology on $J$, and if $J$ is finite, every Grothendieck topology
is of this form (to recover $J'$ from a topology $T$, one takes the
subset $J' \in J$ of elements $j \in J$ such that $T(j) \subset
\Omega(j)$ consists of the maximal sieve $\Y(i)$). We have
$\Shv(J,\E) \cong \Fun({J'}^o,\E)$ for any finitely presentable $\E$,
and the associated sheaf functor $a$ is given by restriction to $J'
\subset J$. If $J$ is left-pointed --- that is, has the smallest
element $o \in J$ --- then the topology corresponding to $\{o\}
\subset J$ is the submaximal topology of Example~\ref{T.exa}.
\end{exa}

\begin{remark}\label{base.rem}
One does not need to specify the whole topology $T \subset \Omega$
to define sheaves and compute the associated sheaf functor. Indeed,
say that a {\em base} $C$ of a topology $T$ on a small category $I$
is a subfunctor $C \subset T$ such that for any $i \in I$, $C(i)^o
\subset T(i)^o$ is a cofinal partially ordered subset. Then for any
functor $E:I^o \to \E$ to some finitely presentable $\E$, one can
replace the colimit over $T(i)^o$ in \eqref{a.eq} with the colimit
over $C(i)^o$, so that a map $E \to a(E)$ is an isomorphism --- that
is, $E$ is a sheaf --- if and only if $\Hom(s,E) \cong E(i)$ for any
$i \in I$ and $s \in C(i)$. Note that $T(i) \subset \Omega(i)$ can
be recovered as the {\em right closure} of $C(i) \subset \Omega(i)$,
that is, the subset of sieves $s \in \Omega(i)$ with non-empty
$C(i)/s$.
\end{remark}

\begin{remark}
For any small $I$, the functor $\Omega:I^o \to \Sets$ has the
following universal property: for any monomorphism $X \to Y$ in
$\Fun(I^o,\Sets)$, there exists a unique map $Y \to \Omega$ that
fits into a cartesian square
\begin{equation}\label{om.sq}
\begin{CD}
X @>>> \ppt_I \\
@VVV @VV{1}V\\
Y @>>> \Omega,
\end{CD}
\end{equation}
where $1$ is the embedding \eqref{min.t.eq} (for a representable
$Y$, this is simply the definition of $\Omega$). Then in particular,
a subfunctor $T \subset \Omega$ defines a map $j:\Omega \to \Omega$,
and one can show that $T$ is a topology if and only if \thetag{i} $j
\circ 1 = 1$, \thetag{ii} $\wedge \circ (j \times j) = j \circ
\wedge$, and \thetag{iii} $j \circ j = j$, where $\wedge:\Omega
\times \Omega \to \Omega$ corresponds to the square \eqref{om.sq}
for the embedding $(1 \times 1):\ppt_I = \ppt_I \times \ppt_I \to
\Omega \times \Omega$.  The functor $\Omega$ is known as the
``subobject classifier'' in $\Fun(I^o,\Sets)$, and the image
$\Omega_j$ of the idempotent endomorphism $j:\Omega \to \Omega$ is a
subobject classifier in the sheaf category $\Shv(I,\Sets)$. In
general, sheaf categories $\Shv(I,\Sets)$ are known as {\em
  toposes}, and it turns out that it is the existence of a subobject
classifier that characterizes them in the most natural way. This is
the subject of the abstract topos theory, for which we refer the
reader to the wonderful book \cite{J} (although when reading it, it
is worthwhile to remember that many of the proofs are based on those
in \cite{SGA4} even when it is not expicitly indicated). The
treatment of topologies in terms of the subobject classifier is
\cite[Chapter 3.1]{J}.
\end{remark}

\subsection{Coverings.}\label{cov.subs}

A practical way to describe a sieve is by specifying a {\em covering
  family} $\{i_\alpha \to i\}$, indexed by some set of indices
$\alpha$, with the sieve $s(\{i_\alpha \to i\})$ consisting of all
maps $i' \to i$ that factor through one of the maps $i_\alpha \to
i$. A {\em Grothendieck pretopology} on $I$ is given by a bunch of
covering families for every $i \in I$, again satisfying some axioms
(that insure in particular that the resulting collection of sieves
is a base for a Grothendieck topology in the sense of
Remark~\ref{base.rem}).

We will not need the full definition of a pretopology (see
e.g.\ \cite[Chapter 0.3]{J}) but we do need one somewhat degenerate
example when all the covering families consist of one
object. Namely, say that a class $F$ of maps in $I$ is a {\em
  covering class} if it is closed under compositions, contains all
the identity maps, and such that a map $f:i' \to i$ in $F$ admits a
pullback with respect to any map $i'' \to i$ (``maps in $F$ admit
pullbacks''), and the pullback is also in $F$ (``$F$ is stable under
pullbacks''). Then for any covering class $F$, setting $C(i) = \{
s(\{f:i' \to i\})| f \in F \}$ defines a base $C$ of a Grothendieck
topology $T$ on $I$ that we call the {\em $F$-topology}. A functor
$E:I^o \to \E$ to some finitely presentable $\E$ is separable
w.r.t.\ the $F$-topology iff for any $f:i' \to i$ in $F$, $E(i) \to
E(i')$ is injective, and it is a sheaf iff the square
\begin{equation}\label{E.sq}
\begin{CD}
E(i) @>>> E(i')\\
@VVV @VVV\\
E(i') @>>> E(i' \times_i i')
\end{CD}
\end{equation}
is cartesian. Up to an equivalence, the category $C(i)$, $i \in I$
is obtained by considering the full subcategory $I /_F i \subset
I/i$ spanned by maps in $F$, and identifying all the maps between
any two objects (so that in the end, there is at most one map
between any two objects, and the category is canonically equivalent
to a partially ordered set). In particular, if an arrow $f:i' \to i$
splits --- that is, admits an inverse $g:i \to i'$, $f \circ g =
\id$, --- then the corresponding object $s(f) \in C(i)$ is
isomorphic to the maximal sieve.

For any object $i \in I$, denote by $\ev_i:\Fun(I^o,\E) \to \E$ the
evaluation functor sending $E$ to $E(i)$, and note that since $\E$
is complete, $\ev_i$ has a right-adjoint $\Y_i:\E \to \Fun(I^o,\E)$
given by
\begin{equation}\label{yo.i.2.eq}
\Y_i(E)(i') = E(\Hom(i,i')), \qquad i' \in I,
\end{equation}
where the right-hand side is shorthand for ``product of copies of
$E$ numbered by elements in the set $\Hom(i,i')$''. Equivalently, we
have $\Y_i(E) \cong Y^i(E)^o$, where $Y^i(E)$ is as in
\eqref{yo.i.eq}. More generally, consider the projective completion
$\Pro(I)$ of the category $I$. Then since $\E$ has filtered
colimits, any functor $E:I^o \to \E$ uniquely extends to a
continuous functor $\wt{E} = \iota_!:\Pro(I)^o = \Ind(I^o) \to \E$,
and we can define a pair of adjoint functors
\begin{equation}\label{ev.Y.eq}
\ev_{\wt{i}}:\Fun(I^o,\E) \to \E, \qquad \Y_{\wt{i}}:\E \to
\Fun(I^o,\E)
\end{equation}
for any proobject $\wt{i} \in \Pro(I)$ by
\begin{equation}\label{ev.Y}
\begin{aligned}
\ev_{\wt{i}}(E) &= \wt{E}(\wt{i}) = \colim_{\wt{i} \to i}E(i),
\qquad E \in \Fun(I^o,\E),\\
\Y_{\wt{i}}(E')(i') &= E'(\Hom(\wt{i},i')) = \lim_{\wt{i} \to
  i}E'(\Hom(i,i')), \qquad E' \in \E,
\end{aligned}
\end{equation}
where the colimit resp.\ the limit actually reduce to a cofinal
filtered diagram representing $\wt{i}$ resp.\ its opposite.

\begin{defn}
A proobject $\wt{i} \in \wc{I}$ is {\em $F$-liftable} if
$\Hom(\wt{i},-)$ sends maps in $F$ to surjective maps.
\end{defn}

\begin{lemma}\label{Y.le}
For any $F$-liftable proobject $\wt{i} \in \wc{I}$, the functor
$\Y_{\wt{i}}$ of \eqref{ev.Y.eq} takes values in $\Shv(I,\E) \subset
\Fun(I^o,\E)$, and the adjoint evaluation functor $\ev_{\wt{i}}$
factors through the associated sheaf functor $a:\Fun(I^o,\E) \to
\Shv(I,\E)$.
\end{lemma}

\proof{} For the second claim, it suffices to check that
$\ev_{\wt{i}}$ inverts the map $E \to a_0(E)$ for any $E$. By
\eqref{ev.Y} and \eqref{a.eq}, we have
\begin{equation}\label{ev.a}
\ev_{\wt{i}}(a_0(E)) \cong \colim_{\wt{i} \to i}\colim_{s \in
  C(i)}\Hom(s,E) = \colim_{(\wt{i} \setminus C)^o}E(s),
\end{equation}
where $\wt{i} \setminus C = (\wt{i} \setminus I) \times_I C$ is the
category of triples $\langle s,i,\wt{i} \to i\rangle$ of an object
$i \in I$, an arrow $\wt{i} \to i$, and a sieve $s \in
C(i)$. However, the $F$-liftability of $\wt{i}$ insures that for any
such $\langle s,i,\wt{i} \to i \rangle$, with $s$ represented by an
arrow $f:i' \to i$ in $F$, the arrow $\wt{i} \to i$ factors through
$f$, and then $f^*(s)$ is split. Therefore the full subcategory in
$(\wt{i} \setminus C)^o$ spanned by triples with split $s$ is
cofinal, so that we can reduce the colimit in \eqref{ev.a} to a
colimit over this subcategory. But the smaller colimit is exactly
$\ev_{\wt{i}}(E)$, so we have proved the second claim. By
adjunction, this means that for any $E' \in \Fun(I^o,\E)$, $E \in
\E$, any map $g:E' \to \Y_{\wt{i}}(E)$ factors uniquely through the
adjunction map $E' \to e(a(E'))$. Taking $E' = \Y_{\wt{i}}(E)$ and
$g=\id$, we deduce that the adjunction map $E' \to e(a(E'))$ splits,
so that $E'$ is a sheaf.
\endproof

\begin{remark}
Note that since the colimit in \eqref{ev.Y} reduces to a filtered
colimit, the functor $\ev_{\wt{i}}$ commutes with finite limits. If
$\E=\Sets$, then this by definition means that the adjoint pair
$\langle \ev_{\wt{i}},\Y_{\wt{i}} \rangle$ then defines a {\em
  point} of the topos $\Shv(I,\Sets)$.
\end{remark}

To construct $F$-liftable proobjects in $I$, take an object $i \in
I$, and let $\Cov(i) \subset I/i$ be the full subcategory spanned by
arrows $i' \to i$ in $F$, with the induced projection
$\sigma(i):\Cov(i) \to I/i \to I$. The class $\sigma(i)^*F$ of maps
$f$ such that $\sigma(i)(f) \in F$ is then a covering class of maps
in $\Cov(i)$, and we have the following.

\begin{defn}
An {\em $F$-hull} of an object $i \in I$ is a
$\sigma(i)^*F$-liftable proobject $\wt{i}$ in $\Cov(i)$.
\end{defn}

\begin{lemma}\label{hull.le}
\begin{enumerate}
\item For any $F$-hull $\wt{i} \in \Pro(\Cov(i))$ of an object $i
  \in I$, the proobject $\sigma(i)(\wt{i}) \in \Pro(I)$ is
  $F$-liftable.
\item For any $i \in I$, there exists an $F$-hull $\wt{i} \in
  \Pro(\Cov(i))$.
\end{enumerate}
\end{lemma}

\proof{} For \thetag{i}, assume given a map $f:i'' \to i'$ in $F$;
we need to show that any map $g:\sigma(i)(\wt{i}) \to i'$ factors
through $i''$. But by definition, $\wt{i}$ comes from a projective
system in $\Cov(i)$, so $g$ factors through $\sigma(g')$ for some
map $g':\wt{i} \to i_0$ in $\Cov(i)$. Then it suffices to check that
$\sigma(g')$ factors through the map $i'' \times_{i'} i_0 \to i_0$;
but this map is in $F$, so it is comes from map in $\Cov(i)$, and
$\wt{i}$ is $\sigma(i)^*F$-liftable in $\Cov(i)$.

For \thetag{ii}, the argument is completely standard and goes back
at least to \cite{toho} (where it is already called
``standard''). Note that since $I$ is small, so is $\Cov(i)$, and
then for any proobject $\wt{i} \in \Pro(\Cov(i))$ there is a set $S$
whose elements $s \in S$ enumerate, up to an isomorphism, all
diagrams
\begin{equation}\label{lft.dia}
\begin{CD}
\wt{i} @>>> i_s' @<{f}<< i_s''
\end{CD}
\end{equation}
in $\Cov(i)$, $f \in \sigma(i)^*F$. Moreover, $\Cov(i)$ has finite
products, and finite products of maps in $F$ are in $F$, so for any
finite subset $S_0 \subset S$, we can define a proobject
$\wt{i}[S_0]$ as the fibered product
$$
\begin{CD}
\wt{i}[S_0] @>>> \prod_{s \in S_0}i_s''\\
@VVV @VVV\\
\wt{i} @>>> \prod_{s \in S_0} i_s'.
\end{CD}
$$
We can then let $H(\wt{i}) = \lim_{S_0 \subset S}\wt{i}[S_0]$, where
the limit is taken over the directed partially ordered set of finite
subsets in $S$. This is a proobject that comes equipped with a map
$H(\wt{i}) \to \wt{i}$, and by construction, for any diagram
\eqref{lft.dia}, the composition map $H(\wt{i}) \to i_s'$ factors through
$f_s$. To finish the proof, let $i_0=i$, with the identity map $i
\to i$, define inductively $i_{n+1} = H(i_n)$, $n \geq 0$, and take
$\wt{i} = \lim_ni_n$.
\endproof

\section{Hypercoverings.}\label{hc.sec}

\subsection{Partially ordered sets.}\label{pos.subs}

Fix a category $I$ equipped with a covering class $F$. Recall that
$\ppt^< \cong [0]^< \cong [1]$ is the single arrow category, and
$[1]^o \cong [1]$, so that coverings of some object $i \in I$ can be
understood as $i$-coaugmented functors from the point category
$\ppt$ to $I$. The notion of a {\em hypercovering} extends this to
categories other than $\ppt$. The usual application is to the
simplex category $\Delta$, see below in
Subsection~\ref{del.subs}, but it is instructive to develop the
theory in a more general context. We start with finite partially
ordered sets.

\begin{defn}\label{hc.pos.def}
For any left-pointed finite partially ordered set $J$, a functor
$E:J^o \to I$ is a {\em hypercovering} if for any non-empty
left-closed subset $J_0 \subset J$, $\lim_{J_0^o} E$ exists, and for
any two non-empty left-closed subsets $J_0 \subset J_1 \subset J$,
the natural map $\lim_{J_1^o}E \to \lim_{J_0^o}E$ is in $F$.
\end{defn}

In the situation of Definition~\ref{hc.pos.def}, we denote by
$\HCov(J) \subset \Fun(J^o,I)$ the full subcategory spanned by
hypercoverings. If $\phi:J' \to J$ is a left-pointed functor between
left-pointed finite partially ordered sets, then \eqref{kan.eq}
immediately shows that for any $E \in \HCov(J')$ the right Kan
extension $\phi^o_*E$ exists and is a hypercovering, and if $\phi$
is a left-closed embedding, then $\phi^{o*}$ also taulotogically
sends hypercoverings to hypercoverings. For any finite partially
ordered set $J$, we have the left-pointed partially ordered set
$J^<$, and we define a {\em $J$-hypercovering} of an object $i \in
I$ as an $i$-coaugmented functor $E:J^{o>} = J^{<o} \to I$ that is a
hypercovering in the sense of Definition~\ref{hc.pos.def}. We denote
the category of $J$-hypercoverings of an object $i$ by $\HCov(J,i)
\subset \HCov(J^<)$. Note that non-empty left-closed subsets in
$J^<$ correspond bijectively to all left-closed subsets in $J$. If
$J=\ppt$, then the only left-closed subsets in $J$ are the empty set
and $J$ itself, so that a $J$-hypercovering is a covering and
$\HCov(\ppt,i) \cong \Cov(i)$. In the general case, for any $j \in
J$, $J/j \subset J$ is left-closed, and since it has the largest
element, \eqref{cof.eq} provides an isomorphism
\begin{equation}\label{F.ev.eq}
\lim_{(J/j)^{<o}}E \cong E(j)
\end{equation}
for any $E \in \HCov(J,i)$. Thus for any $j \leq j'$, the map $E(j')
\to E(j)$ is in $F$, so that a $J$-hypercovering factors through the
subcategory $I_F \subset I$ with the same objects as $I$ and maps that are
in $F$.

\begin{exa}\label{n.hc.exa}
For any $n \geq 0$, any non-empty left-closed subset $J_0 \subset
[n]$ is of the form $[m] = \{0,\dots,m\} \subset [n]$, $0 \leq m
\leq n$, so that $J_0^o$ has the smallest element, and
$\lim_{J_0^{<o}}$ reduces to evaluation by \eqref{F.ev.eq}. Thus
$E:[n]^o \to I$ is a hypercovering if and only if it factors through
$I_F$.
\end{exa}

\begin{remark}\label{hc.12.rem}
A non-empty left-closed subset $J_0 \subset J$ in a left-pointed
partially ordered set $J$ is the preimage of $\{0\} \subset [1]$
under a unique left-pointed map $J \to [1]$, and similarly, pairs
$J_0 \subset J_1 \subset J$ of non-empty left-closed subsets
correspond to left-pointed maps $J \to [2]$. Then by virtue of
\eqref{kan.eq} and Example~\ref{n.hc.exa},
Definition~\ref{hc.pos.def} can be rephrased as follows: $E:J^o \to
I$ is a hypercovering if and only if for any left-pointed map
$\phi:J \to [2]$, $\phi^o_*E$ exists, is universal in the sense of
Remark~\ref{univ.kan.rem}, and factors through $I_F$.
\end{remark}

\begin{exa}\label{v.hc.exa}
Let $J = \{0,\dots,n\}$ be the set of integers $0,\dots,n$ with
discrete order. Then any subset $J_0 \subset J$ is left-closed, and
if $J_0$ has at most one element, $\lim_{J_0^{<o}}$ again reduces to
evaluation. For larger $J_0$ this is not true. However, since maps
in a covering class admits pullbacks and are stable under pullbacks,
induction on the cardinality of $J_0$ again shows that any functor
$E:J^{<o} \to I_F \subset I$ is a hypercovering.
\end{exa}

\begin{exa}
Consider the square $[1]^2$ of the single arrow category $[1]$. Then
$[1]^{2o} \cong [1]^2$, functors $[1]^{2o} \to I$ correspond to
commutative squares
\begin{equation}\label{i.sq}
\begin{CD}
i_{11} @>>> i_{10}\\
@VVV @VVV\\
i_{01} @>>> i_{00},
\end{CD}
\end{equation}
and such a square is a hypercovering iff the maps $i_{01},i_{10} \to
i_{00}$ is in $F$, and so is the map $i_{11} \to i_{10}
\times_{i_{00}} i_{01}$. In particular, it is not sufficient to
require that the functor factors through $I_F$.
\end{exa}

In the general case, testing whether a given functor $E:J^{o>} \to
I$ is a $J$-hypercovering requires checking a lot of
conditions. However, this can be reduced to one condition for each
element $j \in J$ by induction and the following result.

\begin{lemma}\label{hc.le}
Let $j \in J$ be a maximal element in a finite partially ordered set
$J$, let $J' = J \setminus \{j\} \subset J$, and let $L(j) = J/j
\cap J' \subset J$. Then $E:J^{o>} \to I$ is a $J$-hypercovering if
and only if
\begin{enumerate}
\item its restriction to $J' \subset J$ is a $J'$-hypercovering, and
\item the natural map $E(j) = \lim_{(J^</j)^o}E \to
  \lim_{L(j)^{<o}}E$ is in $F$.
\end{enumerate}
\end{lemma}

\proof{} The ``only if'' part is clear. For the ``if'' part, by
induction on cardinality, it suffices to check that $\lim_{J^{o>}}E$
exists, and the map $\lim_{J^{o>}}E \to \lim_{J_0^{o>}}E$ is in $F$
for any left-closed $J_0 \subset J$. For the first claim, let $\V =
\{0,1\}^<$ be as in Example~\ref{V.exa}, and consider the map
$\phi:J \to \V$ sending $L(j)$ to $o$, $j$ to $1$, and the rest to
$0$. Then all the comma-fibers $J^</_{\phi^<} \V$ either have
cardinality smaller than $J$, or have a largest element, or both, so
that $\phi^{<o}_*E$ exists by induction, and satisfies the
assumptions of the Lemma for $J = \V$, $j=1$. We are thus reduced to
this case that immediately follows from Example~\ref{v.hc.exa}. For
the second claim, consider the map $\phi:J \to \V$ sending $J_0 \cap
(J \setminus \{j\})$ to $o$, $j$ to $1$ and the rest to $0$; then
again, $\phi^{<o}_*E$ exists by induction, and we are again reduced
to the obvious case $J=\V$, $j=1$.
\endproof

\begin{exa}\label{hc.top.exa}
The class $F$ of all monomorphisms is a covering class in any
complete category $\E$. In this case, by Lemma~\ref{hc.le}, a
functor $J^o \to \E$ is a hypercovering in the sense of
Definition~\ref{hc.pos.def} if and only if it is a separable
presheaf for the submaximal topology of Example~\ref{T.exa}. More
generally, for any $I$ and $F$, a functor $E:J^o \to I$ is a
hypercovering if and only if $\Hom(s,E)$ defined as the limit
\eqref{hom.eq} exists for any sieve $s \subset \Y(i)$ in the
submaximal topology, and the map $E(i) \to \Hom(s,E)$ is in $F$.
\end{exa}

It is useful to extend the class $F$ to a covering class in the
hypercovering categories $\HCov(J,i)$. Namely, note that the
category $[n]$ is left-pointed for any $n \geq 0$, so that for any finite
partially ordered set $J$, the product $[n] \times J^<$ is a
left-pointed finite partially ordered set. Recall also that $[n]^o
\cong [n]$.

\begin{defn}\label{F.J.def}
A map $f:E_0 \to E_1$ in $\HCov(J^<)$ is in the class $F(J)$ if the
corresponding functor $\iota(f):[1] \times J^{<o} \cong ([1] \times
J^<)^o \to I$ is a hypercovering.
\end{defn}

\begin{lemma}\label{F.le}
For any finite partially ordered set $J$, the class $F(J)$ is a
covering class in $\HCov(J^<)$, and for any $i \in I$, it restricts
to a covering class in $\HCov(J,i) \subset \HCov(J^<)$. For any map
$\phi:J' \to J^<$ from a finite partially ordered set $J'$, the
right Kan extension functor $\phi^{<o}_*$ sends maps in $F(J')$ to
maps in $F(J)$. Moreover, for any left-closed subset $J' \subset J$
with the embedding functor $\phi:J' \to J$, and $E \in \HCov(J^<)$,
the adjunction map $a:E \to \phi^{<o}_*\phi^{<o*}E$ is in $F(J)$.
\end{lemma}

\proof{} For the first claim, denote $J_n = [n] * J$, $n = 0,1$,
where $- * -$ is the extended product \eqref{st.eq}, and note that
for any composable pair $f$, $f'$ of maps in $F(J)$, the functor
$\iota(f,f')$ satisfies the condition \thetag{ii} of
Lemma~\ref{hc.le} at any $l \times i \in J^{<o}_1 \subset ([2]
\times J^<)^o$. Indeed, for $l=0,1$, this is the condition at $l
\times i \in J_0^{<o}$ for $\iota(f')=(t \times \id)^*\iota(f,f')$
on the nose (note that under the identifications $[1]^o \cong [1]$,
$[2]^o \cong [2]$ we have $s^o = t$ and $t^o=s$). For $l=2$, the
functor $s:[1] \to [2]$ has a right-adjoint $s^\dg:[2] \to [1]$, and
the product $s^\dg \times \id$ restricts to a functor $s^\dg \times
\id:L(2 \times j)^{<o} \to L(1 \times j)^{<o} \subset [1] \times
J^{<o}$ right-adjoint to $s \times \id$. We then have
$$
\begin{aligned}
\lim_{L(2 \times j)^{<o}}\iota(f,f') &\cong \lim_{L(1 \times
  j)^{<o}}s^\dg_*\iota(f,f') \cong\\
&\cong \lim_{L(1 \times
  j)^{<o}}s^*\iota(f,f') \cong \lim_{L(1 \times
  j)^{<o}}\iota(f),
\end{aligned}
$$
so we are reduced to Lemma~\ref{hc.le}~\thetag{ii} for $\iota(f)$.
Thus $\iota(f,f')$ is a $J_1$-hyper\-covering, and since
$m:[1] \to [2]$ also has a right-adjoint $m^\dg$, $\iota(f \circ f')
\cong m^*\iota(f,f') \cong m^\dg_*\iota(f,f')$ is a
$J_0$-hypercovering. This shows that $F(J)$ is closed under
compositions. The identification \eqref{F.ev.eq} then immediately
implies that maps in $F(J)$ are pointwise in $F$, so that they admit
pullbacks, and they are stable under pullbacks again by
Lemma~\ref{hc.le}. For the second claim, it suffices to observe that
$\iota(\phi^{<o}_*(f)) \cong (\id \times \phi)^{<o}_*\iota(f)$, and
for the third claim, let $J'' = J'_0 \cup (\{1\} \times J) \subset
J_0$, with the embedding map $\phi_1:J'' \to J_0$, and note that
$\iota(a) \cong \phi^{<o}_{1*}\phi^{<o*}_1\iota(\id_E)$, where
$\id_E:E \to E$ is the identity map.
\endproof

\subsection{Thin $ML$-categories.}\label{ml.subs}

To go beyond partially ordered sets, recall that a {\em
  factorization system} $\langle L,R \rangle$ on a category $\C$ is
given by two classes of maps $L$, $R$ in $\C$, both closed under
compositions and containing all isomorphisms, such that any map $f:c'
\to c$ in $\C$ factors as
\begin{equation}\label{facto.dia}
\begin{CD}
c' @>{l}>> c'' @>{r}>> c,
\end{CD}
\end{equation}
with $l \in L$, $r \in R$, and the factorization is unique up to a
unique isomorphism. This very useful notion goes back to \cite{bou},
and we refer to \cite[Section 2]{bou} for further details. Given a
factorization system $\langle L,R \rangle$, we will denote by
$\C_L,\C_R \subset \C$ the subcategories with the same objects as
$\C$ and maps that are in $L$, $R$, and we note that as a
consequence of the uniqueness of the factorization
\eqref{facto.dia}, the obvious functor $\C_R / c \to \C / c$, $c \in
\C$ is a fully faithful embedding that admits a left adjoint
(sending $f:c' \to c$ to the second term $r:c'' \to c$ in
\eqref{facto.dia}).

\begin{defn}
A {\em thin $ML$-category} is a small category $X$ equipped with a
factorization system $\langle M,L \rangle$ such that for any $x \in
X$, $X_L / x$ is a finite partially ordered set. A functor $\phi:X
\to X'$ between thin $ML$-categories $\langle X,M,L \rangle$,
$\langle X,M',L' \rangle$ is an {\em $ML$-functor} if it sends maps
in $M$ resp.\ $L$ to maps in $M'$ resp.\ $L'$.
\end{defn}

For any thin $ML$-category $X$, we can turn the augmented category
$X^<$ into a thin $ML$-category by putting all the maps $o \to x$,
$x \in X$ in the class $L$, so that $X^<_L \cong (X_L)^<$. For any
$x \in X$, the partially ordered set $X_L / x$ has the largest
element $\id:x \to x$, and we denote $L(x) = (X_L / x) \setminus
\id$.

\begin{defn}\label{hc.ml.def}
For any thin $ML$-category $X$, an {\em $X$-hypercovering} of an
object $i \in I$ is an $i$-coaugmented functor $E:X^{<o} = X^{o>}
\to I$ such that for any $x \in X$, the limit
$\lim_{L(x)^{<o}}\sigma(x)^{<o*}E$ exists, and the map
\begin{equation}\label{hc.ml.eq}
E(x) \cong \lim_{(X^<_L/x)^{<o}}\sigma(x)^{<o*}E \to
\lim_{L(x)^{<o}}\sigma(x)^{<o*}E
\end{equation}
is in the class $F$.
\end{defn}

\begin{exa}\label{pos.ml.exa}
Any finite partially ordered set $J$ is trivially a thin
$ML$-category, with $M$ resp.\ $L$ consisting of the identity maps
resp.\ all maps. In this case, Lemma~\ref{hc.le} immediately shows
that Definition~\ref{hc.ml.def} reduces to
Definition~\ref{hc.pos.def}.
\end{exa}

\begin{remark}\label{hc.top.rem}
One can also describe hypercoverings in terms of Grothen\-dieck
topologies, as in Example~\ref{hc.top.exa}. Namely, for any thin
$ML$-category $X$, a sieve $s$ on an object $x \in X^<_L$ defines an
``induced'' sieve on the same object in $X^<$ consisting of all maps
$x' \to x$ whose component $x'' \to x$ of the factorization
\eqref{facto.dia} is in $s$. Then all induced sieves $s \in \Y(x)$
for the submaximal topology on $X^<_L$ form a topology on $X^<$, and
$E:X^{<o} \to I$ is a hypercovering iff $\Hom(s,E)$ exists for any
$s$, and the map $E(x) \to \Hom(s,E)$ is in $F$. In good cases, this
induced topology can be described directly; for example, one can
show that if all maps $l \in L$ are monomorphisms, and any map $m
\in X$ admits a one-sided inverse $l$, $m \circ l = \id$, then the
induced topology on $X^<$ is submaximal (that is, consists of all
non-empty sieves). This happens for the category $\Delta$ considered
below in Subsection~\ref{del.subs}.
\end{remark}

For an arbitrary thin $ML$-category $X$, we denote the category of
$X$-hypercoverings of an object $i \in I$ by $\HCov(X,i)$. By virtue
of Example~\ref{pos.ml.exa}, this is consistent with our earlier
notation. For any map $l:x' \to x$ in the class $L$, we have
$(X_L/x)/l \cong X_L/x'$, so that $\sigma(x)^{<o*}$ sends
$X$-hypercoverings to $(X_L/x)$-hypercoverings. In particular, this
implies that any $X$-hypercovering $E:X^{<o} \to I$ sends $X^{<o}_L
\subset X^{<o}$ into $I_F \subset I$, so that we have a full
embedding $\HCov(X,i) \subset \Fun(X^o,\Cov(i))$.

\begin{exa}
Say that a full subcategory $X \subset X'$ in a thin $ML$-category
$\langle X,M,L \rangle$ is an {\em $ML$-subcategory} if for any map
$f:x' \to x$ in $X'$, the middle term $x''$ of its decomposition
\eqref{facto.dia} is in $X'$. Then $X'$ with the classes $M$, $L$ is
a thin $ML$-category, and the embedding $X' \to X$ is an
$ML$-functor.
\end{exa}

\begin{exa}\label{prod.ml.exa}
For any two thin $ML$-categories $\langle X,M,L \rangle$, $\langle
X',M',L' \rangle$, the product $\langle X \times X',M \times M',L
\times L'\rangle$ is a thin $ML$-category, and so is the extented
product $X * X'$ of \eqref{st.eq}.
\end{exa}

Combining Example~\ref{pos.ml.exa} and Example~\ref{prod.ml.exa}, we
see that the extended product $J * X$ of a finite partially
ordered set $J$ and a thin $ML$-category $X$ is naturally a thin
$ML$-category. In particular, as in Definition~\ref{F.J.def}, we can
take $J=[0]$, and say that a map $f$ in $\HCov(X,i)$ is in the class
$F(X)$ if $\iota(f)$ is a hypercovering. Then Lemma~\ref{F.le}
immediately shows that $F(X)$ is a covering class.

\begin{exa}\label{Pos.ml.exa}
Let $\Pos$ be the category of all non-empty finite partially ordered
sets, let $M$ be the class of all surjective maps, and let $L$ be
the class of all injective maps $J' \to J$ that are full when
considered as functors (that is, $J'$ is identified with its image
in $J$ equipped with the induced partial order). Then $\Pos$ is a
thin $ML$-category in the sense of Definition~\ref{hc.ml.def}.
\end{exa}

\begin{lemma}\label{ev.le}
Assume given an object $x \in X$ in a thin $\Hom$-finite
$ML$-category $\langle X,M,L \rangle$, and let $\phi:\ppt \to X$ be
the embedding onto $x$. Then for any object $i \in I$, $\phi^{<o*}$
and $\phi^{<o}_*$ define an adjoint pair of functors between
$\HCov(X,i)$ and $\Cov(i)$ sending maps in $F(X)$, $F$ to maps in
$F$, $F(X)$.
\end{lemma}

\proof{} Since the restriction of an $X$-hypercovering to $X^{<o}_L$
factors through $I_F$, $\phi^{<o*}$ sends hypercoverings to
coverings. Moreover, by the definition of a covering class,
$\Cov(i)$ has finite products. Since $X$ is $\Hom$-finite, the right
Kan extension $\phi^{<o}_*E$ then exists by \eqref{kan.eq} for any
$E \in \Cov(i)$ and is given by
$$
\phi^{<o}_*E(x') \cong \prod_{f:x \to x'}E,
$$
the product of copies of $E$ numbered by maps $f:x \to x'$. But the
set of maps $x \to x'$ splits into a disjoint union according to the
isomorphism class of the middle term $x''$ of the decomposition
\eqref{facto.dia}, and $\lim_{L(x')^{<o}}\phi^{<o}_*E$ is then given
by the same product but over the subset of maps such that $x'' \to
x'$ is not a isomorphism. Since maps in $F$ are stable under
pullbacks, the map \eqref{hc.ml.eq} for $\phi^{<o}_*E$ is in $F$ for
any $x' \in X$.
\endproof

\begin{lemma}\label{M.le}
Assume given a thin $ML$-category $\langle X,M,L \rangle$ and a full
$ML$-subcategory $X' \subset X$. Let $\phi:X' \to X$ be the
embedding functor, and assume that $X'_M \subset X_M$ is
right-closed. Then for any $X'$-hypercovering $E$, $\phi^{<o}_*E$
exists and is an $X$-hypercovering.
\end{lemma}

\proof{} For any $x \in X$, let $\phi_x:X'_L/x \to X_L/x$ be the
embedding induced by $\phi$, and consider the base change map
\begin{equation}\label{bc.ml.eq}
\sigma(x)^{<o*}\phi^{<o}_*E \to \phi_{x*}^{<o}\sigma(x)^{<o*}E.
\end{equation}
Since $\sigma(x)^{<o*}E$ is a hypercovering, we already know that
the target of \eqref{bc.ml.eq} exists and is a hypercovering, so it
suffices to show that the map is an isomorphism. By \eqref{kan.eq},
this amounts to checking that for any $x \in X$,
${X'}^{<o}_L/_{\phi^{<o}} x$ is cofinal in ${X'}^{<o}/_{\phi^{<o}}
x$. But since $X'_M \subset X_M$ is right-closed, the functor
$X^{<o}/x \to X^{<o}_L/x$ left-adjoint to the embedding $X^{<o}_L/x
\subset X^{<o}/x$ sends ${X'}^{<o}/_{\phi^{<o}} x$ into
${X'}^{<o}_L/_{\phi^{<o}} x$, so we are done by \eqref{cof.eq}.
\endproof

\begin{lemma}\label{adj.ml.le}
Assume given thin $ML$-categories $\langle X,M,L \rangle$, $\langle
X',M',L' \rangle$, and an $ML$-functor $\phi:X \to X'$ that admits a
left-adjoint $\psi:X' \to X$. Moreover, assume that $\psi$ is also
an $ML$-functor. Then $\phi^{<o*}$ sends $X'$-hypercoverings to
$X$-hypercoverings.
\end{lemma}

\proof{} For any $x \in X$, the $ML$-functor $\phi$ induces a
functor $\phi_x:X_L/x \to X'_{L'}/\phi(x)$, and since its adjoint
$\psi$ is also an $ML$-functor, it induces a functor
$\psi_x:X'_{L'}/\phi(x) \to X_L/\psi(\phi(x)) \to X_L/x$
left-adjoint to $\phi_x$. Then $\phi_x^{<o*} \cong \psi_{x*}^{<o}$
sends hypercoverings to hypercoverings for any $x$, and then
so does $\phi^{<o*}$.
\endproof

\begin{lemma}\label{prod.ml.le}
For any two thin $ML$-categories $X_0$, $X_1$, a functor $X_0^{<o}
\times X_1^{<o} \to I$ is a hypercovering if and only if the corresponding
functor $X_0^{<o} \to \Fun(X_1^{<o},I)$ factors through $\HCov(X_1)$
and is a hypercovering with respect to the class $F(X_1)$.
\end{lemma}

\proof{} Since for any $x_0 \times x_1 \in X_0^< \times X^<_1$ we
have $(X_0^< \times X_1^<)_L/(x_0 \times x_1) \cong (X_{0L}^< / x_0)
\times (X_{1L}^< / x_1)$, it suffices to prove the claim when $X_0$
and $X_1$ are finite partially ordered sets. Then
Remark~\ref{hc.12.rem} and \eqref{kan.eq} immediately reduce us to
the case $X_0 = [1]$, and the argument in this case is the same as
in the proof of Lemma~\ref{F.le}.
\endproof

\subsection{Simplicial objects.}\label{del.subs}

Now as usual, denote by $\Delta \subset \Pos$ the full subcategory
spanned by ordinals $[n]$, $n \geq 0$. A {\em simplicial object} in
a category $I$ is by definition a functor $i_\idot:\Delta^o \to I$,
with $i_n = i_\idot([n])$, $n \geq 0$. One traditionally denotes
$\Delta^o I = \Fun(\Delta^o,I)$. An $i$-augmented simplicial object,
$i \in I$ is an $i$-coaugmented functor $i_\idot:\Delta^{<o} \cong
\Delta^{o>} \to I$, and it is explicitly given by a triple $\langle
i_\idot,i,\alpha\rangle$, where $i_\idot$ is a simplicial object,
and $a:i_\idot \to i$ is the augmentation map to the constant
simplicial object with value $i$.

The full subcategory $\Delta \subset \Pos$ inherits a structure of a
thin $ML$-category of Example~\ref{Pos.ml.exa}, so it makes sense to
speak of $\Delta$-hypercoverings in the sense of
Definition~\ref{hc.ml.def}. These are usually just called {\em
  hypercoverings}, and appear in the literature in many places and
forms. Here are some of them.
\begin{enumerate}
\item For any $n \geq 0$, let $\Delta^<_{< n} \subset \Delta^<$ be the
  full subcategory spanned by $[m]$ with $m < n$, and let
  $j_n:\Delta^<_{< n} \to \Delta^<$ be the embedding functor. The
  {\em $n$-th coskeleton} $\cosk_n i_\idot$ of an augmented
  simplicial object $i_\idot:\Delta^{<o} \to I$ is given by the
  right Kan extension $j_{n*}^oj_n^{o*}i_\idot$ if it exists, and if
  it does, it comes equipped with the adjunction map $E \to
  \cosk_nE$. The $0$-th coskeleton $\cosk_0i_\idot$ always exists
  --- by \eqref{kan.eq}, this is just the constant functor with
  value $i = i_\idot(o)$, and the adjunction map $i_\idot \to
  \cosk_0i_\idot$ is $\id$ over $o$ and the augmentation map
  $a:i_\idot \to i$ over $\Delta^o$. Then Lemma~\ref{M.le} applies
  to the embeddings $\Delta_{<n}^< \subset \Delta^<$ and shows that
  an augmented simplicial object $i_\idot:\Delta^{<o} \to I$ is a
  $\Delta$-hypercovering if and only if $\cosk_ni_\idot$ exists for
  any $n \geq 0$, and the map $i_n \to (\cosk_ni_\idot)_n$ is in
  $F$. If $I$ is the category of schemes, and $F$ is the covering
  class of proper maps, then these are the original hypercoverings
  of \cite{Del}.
\item Alternatively, one can describe $\Delta$-hypercoverings in
  terms of sieves, as in Remark~\ref{hc.top.rem}. Namely, let
  $\Delta_n = \Hom(-,[n]):\Delta^o \to \Sets$ be the elementary
  simplex, let $\SS_{n-1} \subset \Delta_n$ be the standard
  simplicial sphere, and augment both by the one-point set
  $\ppt$. Then $i_\idot:\Delta^{<o} \to I$ is a hypercovering if and
  only if for any $n \geq 1$, $\Hom(\SS_{n-1},i_\idot)$ exists and
  the map $i_n = \Hom(\Delta_n,i_\idot) \to \Hom(\SS_{n-1},i_\idot)$
  is in $F$. If this holds, then the same holds for any non-empty
  sieve $\SS \subset \Delta_n$. For example, if $I=\Sets$ is the
  category of sets, and $F$ is the class of all surjective maps,
  then an augmented simplicial set $X_\idot = \langle X_\idot,X,a
  \rangle:\Delta^{<o} \to \Sets$ is a $\Delta$-hypercovering if and
  only if $a:X_\idot \to X$ is a trivial fibration with respect to
  the Kan-Quillen model structure on $\Delta^o\Sets$. Explicitly,
  the hypercovering condition at $[n] \in \Delta$ says that any
  augmented map $\SS_n \to X_\idot$ extends to an augmented map
  $\Delta_n \to X_\idot$. If this holds, then the same is true for
  any non-empty sieve $\SS \subset \Delta_n$. We also note that the
  coaugmentation given by a $\Delta$-hypercovering
  $X_\idot:\Delta^{<o>} \to \Sets$ is universal: if we denote
\begin{equation}\label{pi.0.eq}
\pi_0(X_\idot) = \colim_{\Delta^o}X_\idot, \qquad X_\idot \in
\Delta^o\Sets,
\end{equation}
then for any $\Delta$-hypercovering $\langle X_\idot,X,a \rangle$,
we have $X \cong \pi_0(X_\idot)$, and $a:X_\idot \to X$ is the
canonical map.
\item Now assume that $I$ is a model category, and $F$ is the class
  of fibrations resp.\ trivial fibrations. Then $\Delta$ has the
  structure of a {\em Reedy category}, \cite[Section 5.2]{hov}, with
  our classes $M$ and $L$ corresponding to matching and latching
  maps (and this is why we use this notation). The category of
  simplicial objects in $I$ carries a Reedy model structure, and an
  augmented object $\langle i_\idot,i,a \rangle$ in $I$ is a
  $\Delta$-hypercovering if and only if $a:i_\idot \to i$ is a
  fibration resp.\ trivial fibration. The partially ordered set
  $L([n])$ of Definition~\ref{hc.ml.def} is the latching category
  for the Reedy structure on $\Delta$ (that becomes the matching
  category for the opposite category $\Delta^o$).
\end{enumerate}
From now on, we will also shorten ``$\Delta$-hypercovering'' to
``hypercovering'', and denote $\HCov(i) = \HCov(\Delta,i)$. We
recall that $\HCov(i)$ comes equipped with a covering class
$F(\Delta)$. Recall also that for any two simplicial objects $E_0,E_1
\in \Fun(\Delta^o,\E)$ in a category $\E$, we also have the
simplicial set $\Hhom(E_0,E_1)$ given by \eqref{hhom.eq}, and since
$[0] \in \Delta^o$ is the initial object, we have $\Hom(E_0,E_1)
\cong \lim_{\Delta^o}\Hhom(E_0,E_1) \cong \Hhom(E_0,E_1)_0$. In
particular, this applies to hypercoverings $i_\idot \in \HCov(i)
\subset \Fun(\Delta^o,\Cov(i))$.

\begin{defn}\label{hc.def}
For any object $i \in I$, objects in the category $HC(i)$ are
hypercoverings $i_\idot \in \HCov(i)$, and morphisms are given by
$$
\Hom(i_\idot,i'_\idot) = \pi_0(\Hhom(i_\idot,i'_\idot)),
$$
where $\pi_0(-)$ is as in \eqref{pi.0.eq}.
\end{defn}

To see the set $\pi_0(X_\idot)$ more explicitly, note that since any
object $[n] \in \Delta$ admits a map $[0] \to [n]$, the natural map
$X_0 \to \pi_0(X)$ is surjective. Then define an {\em elementary
  homotopy} between two elements $x,x' \in X_0$ as an element
$\wt{x} \in X_1$ with $X(s)(\wt{x})=x$, $X(t)(\wt{x}) = x'$, and say
that two elements $x,x' \in X_0$ are {\em chain-homotopic} if they
can be connected by a finite chain of elementary homotopies (going
in either direction). Then being chain-homotopic is an equivalence
relation on $X_0$, and $\pi_0(X)$ is the corresponding set of
equivalence classes.

\begin{defn}\label{del.ho.def}
For any two simplicial objects $E_0,E_1 \in \Fun(\Delta^o,\E)$ in a
category $\E$, two maps $E_0 \to E_1$ are {\em chain-homotopic} if
so are the corresponding elements in $\Hhom(E_0,E_1)_0$.
\end{defn}

In this language, morphisms in the category $HC(i)$ of
Definition~\ref{hc.def} are maps in $\HCov(i)$ considered up to a
chain homotopy.

\begin{lemma}\label{copr.le}
Let $\delta^<:\Delta^< \to \Delta^< \times \Delta^<$ be the diagonal
embedding. Then for any hypercovering $E$, $\delta^{<o}_*E$ exists
and is a $(\Delta * \Delta)$-hypercovering.
\end{lemma}

\proof{} The embedding $\phi:\Delta \to \Pos$ satisfies the
assumptions of Lemma~\ref{M.le}, and so does the embedding $\phi
\times \phi:\Delta \times \Delta \to \Pos \times \Pos$. Therefore it
suffices to prove the statement for the diagonal embedding
$\delta:\Pos^< \to \Pos^< \times \Pos^<$, and then apply it to
$\phi^{<o}_*E$ that exists by Lemma~\ref{M.le}. But if we
reinterpret $\Pos^<$ as the category of all finite partially ordered
sets by treating the initial object $o$ as the empty set, then
$\delta$ has a right-adjoint $\mu:\Pos^< \times \Pos^< \to \Pos^<$
given by the cartesian product, and we are done by
Lemma~\ref{adj.ml.le}.
\endproof

\begin{corr}\label{filt.corr}
For any $i \in I$, the category $HC(i)^o$ opposite to $HC(i)$ of
Definition~\ref{hc.def} is filtered.
\end{corr}

\proof{} Definition~\ref{filt.def}~\thetag{i} is obvious ($\HCov(i)$
has finite products). For \thetag{ii}, note that for any
$i_\idot,i'_\idot \in \HCov(i)$ and $n \geq 0$, we have
$$
\Hhom(i'_\idot,i_\idot)_n \cong
\Hom(i'_\idot,\Hhom(\Delta_n,i_\idot)) \cong
\Hom(i'_\idot,\eps_n^*\delta_*i'_\idot),
$$
where $\Hhom(\Delta_n,-)$ is as in \eqref{hhom.E}, and
$\eps_n:\Delta^< \to \Delta^< \times \Delta^<$ is the embedding onto
$[n] \times \Delta^<$. Then by Lemma~\ref{copr.le},
Lemma~\ref{prod.ml.le} and Lemma~\ref{ev.le},
the object $\Hhom(\Delta_n,i_\idot)$ is a hypercovering, and the map
$\Hhom(\Delta_n,i_\idot) \to \Hhom(\SS_{n-1},i_\idot)$ is in
$F(\Delta)$. Thus if are given two maps $f,f':i'_\idot \to i_\idot$
in $\HCov(i)$, we can construct the fibered product
$$
\begin{CD}
i''_\idot @>>> \Hhom(\Delta_1,i_\idot)\\
@V{g}VV @VVV\\
i'_\idot @>{f \times f'}>> \Hhom(\SS_0,i_\idot) \cong i_\idot \times_i
i_\idot
\end{CD}
$$
in $\HCov(i)$, and then $f \circ g$, $f' \circ g$ are connected by
an elementary homotopy, thus give the same map in $HC(i)$.
\endproof

\begin{remark}\label{trunc.rem}
Say that a hypercovering $i_\idot \in \HCov(i)$ is {\em
  $n$-truncated} for some $n \geq 0$ if the map $i_\idot \to
\cosk_n(i_\idot)$ is an isomorphism. Then a $1$-truncated
hypercovering of an object $i$ is the same thing as a covering, and
moreover, any two maps between $1$-truncated hypercoverings are
related by an elemetary homotopy. Thus in fact, the full subcategory
$HC(i)_1 \subset HC(i)$ spanned by $1$-truncated hypercoverings is
canonically identified with $C(i)$.
\end{remark}

\section{Dold-Kan equivalence and derived functors.}\label{dold.sec}

\subsection{Dold-Kan equivalence.}

Let now $\E$ be an additive Karoubi-closed category (for example, an
abelian one). Then we have the {\em Dold-Kan equivalence}
\begin{equation}\label{DK.eq}
\Fun(\Delta^o,\E) \cong C_{\geq 0}(E), \qquad E \mapsto
C_\idot(E),
\end{equation}
where $C_{\geq 0}(\E)$ is the category of chain complexes in $\E$
concentrated in non-negative homological degrees, and $C_\idot(E)$
is the normalized chain complex of the simplicial object $E$. The
original reference for \eqref{DK.eq} is \cite{dold} but there are
many expositions in the literature. A constant functor $\Delta^o \to
\E$ with some value $E$ corresponds to the complex that consists of
$E$ placed in degree $0$. If $\E$ is cocomplete, then by adjunction,
$\colim_{\Delta^o} E$ for some $E \in \Fun(\Delta^o,\E)$ is the
cokernel of the differential $C_1(E) \to C_0(E)$ in the
corresponding chain complex. Augmented simplicial objects correspond
to augmented complexes --- that is, triples $\langle
M_\idot,M,\alpha \rangle$ of a chain complex $M_\idot$, an object
$M$, and a map $\alpha:M_0 \to M$ such that $\alpha \circ d=0:M_1
\to M$. For any $n \geq 0$, \eqref{DK.eq} also identifies
$\Fun(\Delta_{< n}^o,\E)$ with chain complexes $C_{[0,n]}(\E)$
concentrated in degrees $0 \leq i < n$, and the restriction
$j_n^{o*}$ sends a complex $M_\idot$ to its $n$-th stupid truncation
--- that is, to the complex $M_n \to \dots \to M_0$. If $\E$ has
kernels, then the right Kan extension $j^o_{n*}$ exists and sends a
complex $M_n \to \dots \to M_0$ to
\begin{equation}\label{trunc}
\begin{CD}
\Ker d @>>> M_n @>{d}>> \dots @>{d}>> M_0.
\end{CD}
\end{equation}
One can also iterate \eqref{DK.eq} and obtain an equivalence
\begin{equation}\label{DDK.eq}
\Fun(\Delta^o \times \Delta^o,\E) \cong \Fun(\Delta^o,C_{\geq
  0}(\E)) \cong C_{\geq 0,\geq 0}(\E)
\end{equation}
whose target is the category of bicomplexes in $\E$ concentrated in
non-negati\-ve homological bidegrees. By abuse of notation, for any
bicomplex $M_{\idot,\idot}$ corresponding to a bisimplicial object
$M$ under \eqref{DDK.eq}, we will denote by
$\delta^*M_{\idot,\idot}$ the complex corresponding to $\delta^*M$,
where $\delta:\Delta^o \to \Delta^o \times \Delta^o$ is the diagonal
embedding. The complex $\delta^*M_{\idot,\idot}$ is different from
the totalization $\Tot(M_{\idot,\idot})$ of the bicomplex
$M_{\idot,\idot}$ but they are canonically quasiisomorphic ---
namely, there are functorial {\em shuffle maps}
\begin{equation}\label{sh.eq}
\Tot(M_{\idot,\idot}) \to \delta^*M_{\idot,\idot} \to
\Tot(M_{\idot,\idot})
\end{equation}
whose composition is the identity map, and both are
quasiisomorphisms if $\E$ is abelian (again, there are many
expositions of this in the literature, e.g.\ see \cite[Section
  3.4]{trace} for a coordinate-free construction). The individual
terms $M_n$ of the complex $M_\idot = \delta^*M_{\idot,\idot}$ are
finite sums of the terms of the bicomplex
$M_{\idot,\idot}$ --- in particular, if all $M_{\idot,\idot}$ are
projective in $\E$, then so are all the $M_\idot$.

\begin{defn}\label{bimod.def}
Assume that $\E$ is abelian. Then a map $M_{\idot,\idot} \to
N_{\idot,\idot}$ in $C_{\geq 0,\geq 0}(\E)$ is a {\em left
  resp.\ right quasiisomorphism} if $M_{n,\idot} \to N_{n,\idot}$
resp.\ $M_{\idot,n} \to N_{\idot,n}$ is a quasiisomorphism for any
$n \geq 0$.
\end{defn}

\begin{exa}\label{bimod.exa}
The totalization functor $\Tot:C^{\geq 0,\geq 0}(\E) \to C^{\geq
  0}(\E)$ admits two obvious one-sided inverses $\Ll,\Rr:C_{\geq
  0}(\E) \to C_{\geq 0,\geq 0}(\E)$ given by $\Ll(M_\idot)_{i,j} =
M_i$ if $j=0$ and $0$ otherwise, and dually, $\Rr(M_\idot)_{j,i} =
M_i$ if $j=0$ and $0$ otherwise (in terms of the Dold-Kan
equivalence, we have $\Ll \cong \pi_0^*$, $\Rr \cong \pi_1^*$, where
$\pi_0,\pi_1:\Delta^o \times \Delta^o \to \Delta^o$ are the
projections onto the two factors). There are no maps between $\Ll$
and $\Rr$. However, the totalization functor $\Tot$ also admits a
left-adjoint $\I:C_{\geq 0}(\E) \to C_{\leq 0,\geq 0}(\E)$ given by
$\I(M_\idot)_{i,j} = M_{i+j} \oplus M_{i+j+1}$, with both
differentials equal to $d + \id$. Then the isomorphisms $\Tot \circ
\Ll \cong \Tot \circ \Rr \cong \id$ provide by adjunction maps $\I
\to \Ll$ resp.\ $\I \to \Rr$, and is $\E$ is abelian, these two
adjoint maps are a left resp.\ a right quasiisomorphism in the sense
of Definition~\ref{bimod.def}.
\end{exa}

\begin{lemma}\label{DK.le}
Assume that $\E$ is abelian, we have bicomplexes
$M_{\idot,\idot}$, $M'_{\idot,\idot}$ in $\E$, and a map
$f:M_{\idot,\idot} \to M'_{\idot,\idot}$ that is a left or right
quasiisomorphism in the sense of Definition~\ref{bimod.def}. Then
$\delta^*(f):\delta^*M_{\idot,\idot} \to \delta^*M'_{\idot,\idot}$
is also a quasiisomorphism.
\end{lemma}

\proof{} The totalization $\Tot(f)$ is obviously a quasiisomorphism,
and so are the shuffle maps \eqref{sh.eq} both for $M_{\idot,\idot}$
and $M'_{\idot,\idot}$.
\endproof

Dually, one can replace $\Delta^o$ in \eqref{DK.eq} with $\Delta$,
and then $C_{\geq 0}(-)$ gets replaced with the category $C^{\geq
  0}(-)$ of complexes concentrated in non-negative cohomological
degrees --- to deduce the corresponding statement for $\E$, one has
to apply \eqref{DK.eq} to $\E^o$ that is additive and
Karoubi-closed. The rest of the material up to and including
Example~\ref{bimod.exa} and Lemma~\ref{DK.le} also has an obvious
dual counterpart.

\subsection{Covering classes.}

Next, we want to equip our additive Karoubi-closed category $\E$
with a covering class $F$. There are two possible choices. Firstly,
assume that $\E$ has kernels. Then $\E$ is finitely complete, and
we can take the class of all maps. The hypercovering condition for
this class is empty --- $\HCov(M)$ is simply the category of all
chain complexes equipped with an augmentation $M_\idot \to M$ ---
but Definition~\ref{hc.def} is still non-trivial. To describe what
it says, note that any additive category $\E$ is trivially a module
category over the monoidal category $\Z\amod^{ff}$ of finite free
$\Z$-modules. While the Dold-Kan equvalence is not tensor, for any
$M:\Delta^o \to \E$, $V:\Delta^o \to \Z\amod^{ff}$ we have shuffle maps
\begin{equation}\label{shuff}
C_\idot(M) \otimes C_\idot(V) \to C_\idot(M \otimes V) \to
C_\idot(M) \otimes C_\idot(V)
\end{equation}
induced by \eqref{sh.eq} (indeed, if one considers the bicomplex
$C_\idot(M) \boxtimes C_\idot(V)$ with terms $C_n(M) \otimes
C_n(V)$, then $C_\idot(M) \otimes C_\idot(V) \cong \Tot(C_\idot(M)
\boxtimes C_\idot(V))$ and $C_\idot(M \otimes V) \cong
\delta^*(C_\idot(M) \boxtimes C_\idot(V))$). If we let
$\Z[\Delta_1]:\Delta^o \to \Z\amod$ be the simplicial $\Z$-module
generated pointwise by the elementary $1$-simplex, then giving an
elementary homotopy of Definition~\ref{hc.def} between two maps
$f,f':E_\idot \to E'_\idot$ is equivalent to giving a map
\begin{equation}\label{h.0}
h:C_\idot(E_\idot \otimes \Z[\Delta_1]) \to C_\idot(E_\idot')
\end{equation}
with given restriction to $C_\idot(E_\idot \otimes \Z[\SS_0]) \cong
C_\idot(E_\idot) \oplus C_\idot(E_\idot)$. However,
$C_\idot(\Z[\Delta_1])$ is the length-$2$ complex $\Z \to \Z \oplus
\Z$, so that giving chain homotopy between $f$ and $f'$ in the usual
sense is equivalent to giving a map
\begin{equation}\label{h.1}
h':C_\idot(E_\idot) \otimes C_\idot(\Z[\Delta_1]) \to C_\idot(E_\idot').
\end{equation}
By virtue of \eqref{shuff}, the source of \eqref{h.1} is a retract
of the source of \eqref{h.0}, so that an elementary homotopy exists
iff the maps are chain-homotopic in the stardard sense. This
relation is already transitive, so that being chain-homotopic in the
usual sense and in the sense of Definition~\ref{hc.def} amounts to
the same thing.

A more interesting alternative is to assume further that $\E^o$ is
either abelian, or small and preabelian in the sense of
Definition~\ref{preab.def}, so that $\Pro(\E) \cong \Ind(\E^o)^o$ is
abelian, and to take as $F$ the class of all epimorphisms in $\E$.
Since epimorphisms in an abelian category are obviously stable under
pullbacks, and the full embedding $\E \subset \Pro(\E) =
\Ind(\E^o)^o$ reflects epimorphisms, $F$ is again a covering
class. Then \eqref{trunc} immediately shows that an augmented
complex $\langle M_\idot,M,a \rangle$ is an $F$-hypercovering if and
only if $a:M_\idot \to M$ is a quasiisomorphism --- that is,
$M_\idot$ is a left resolution of $M$ in $\Pro(\E) \supset \E$ ---
and $HC(M)$ is the category of left resolutions of the object $M \in
\E$, and chain-homotopy classes of maps between them. Dually, if
$\E$ is preabelian, then hypercoverings in the opposite category
$\E^o$ correspond to right resolutions in $\Ind(\E^o)$.

If we have categories $I$, $I'$ with covering classes $F$, $F'$, and
a functor $\gamma:I \to I'$ that commutes with finite limits and
sends maps in $F$ to maps in $F'$, then it tautologically sends
$F$-hypercoverings to $F'$-hypercoverings. For example, this happens
if $I=\E$ is an abelian category with some projective object $E$,
$E' = \Sets$, $F$ and $F'$ consist of epimorphisms, and $\gamma$ is
the functor $\Hom(E,-)$. If $\E$ has arbitrary sums, then this
particular functor has a left-adjoint $\Sets \to \E$, $S \mapsto
E[S]$, where as in Example~\ref{shv.ab.exa}, $E[S]$ is shorthand for
``the sum of copies of $E$ numbered by elements $s \in S$''. This
adjoint still sends epimorphisms to epimorphisms, even if $E$ is not
projective, but it certainly never commutes with limits. However,
the following is still true.

\begin{lemma}\label{hur.le}
Assume given a set $S$, a collection $E_s \in \E$, $s \in S$ of
objects in an abelian category $\E$ with arbitrary sums, and a
collection of augmented simplicial sets $X_s:\Delta^{<o} \to \Sets$,
$s \in S$ that are hypercoverings with respect to the class of
epimorphisms. Then the augmented simplicial object
$\bigoplus_sE_s[X_s]:\Delta^{<o} \to \E$ is a hypercovering with
respect to the class of epimorphisms.
\end{lemma}

\proof[Sketch of a proof.] Since $\bigoplus_sE_s[X_s]$ is a retract
of $E[X]$, $E = \bigoplus_sE_s$, $X = \coprod_s X_s$, it suffices to
consider the case when $S = \ppt$ has a single element. Moreover, it
suffices to prove that for any $E' \in \E$, the augmented functor
$\Hom(E(X),E'):\Delta^{<o} \to \Ab^o$ is a hypercovering in the
category $\Ab^o$ opposite to the category of abelian groups, again
with respect to epimorphisms, and this amounts to checking the
statement after replacing $\E$ by $\Ab$ and $E$ by $\Hom(E,E')$. In
other words, we might assume right away that $\E=\Ab$. Then $E[X]
\cong E \otimes_{\Z} \Z[X]$, so we may assume that $E=\Z$, and we
may further apply the Dold-Kan equivalence to replace $E[X]$ with
the normalized chain complex $C_\idot(X,\Z)$. The statement is then
well-known: a hypercovering is in particular a weak equivalence, and
$C_\idot(-,\Z)$ sends weak equivalences to quasiisomorphisms. Here
is one argument for this: since all simplicial sets are cofibrant,
any weak equivalence is a composition of trivial cofibrations and
one-sided inverses to such, so it suffices to treat trivial
cofibrations; since filtered colimits in $\Ab$ preserve
quasiisomorphisms, it further suffices to consider elementary
trivial cofibrations, that is, pushouts of horn extensions; but a
pushout of an injective quasiisomorphism is a
quasiisomorphism, so it further suffices to consider horn embeddings.
These give quasiisomorphisms simply by the way the Dold-Kan
equivalence works.
\endproof

\begin{remark}
The statement of Lemma~\ref{hur.le} is essentially a fact from
topology (``no homotopy implies no homology'' -- roughly speaking,
the Hurewicz theorem). The proof sketched above uses Quillen model
structures. There are many alternative proofs, but none are
elementary and/or purely categorical.
\end{remark}

\subsection{Dold-derived functors.}

Now assume given a functor $E:\C \to \E$ between complete categories
$\C$, $\E$, and extend $E$ to a functor $E:\Fun(\Delta^o,\C) \to
\Fun(\Delta^o,\E)$ by applying it termwise. The following is
essentially a classic observation due to Dold.

\begin{lemma}\label{ch.le}
If two maps $f$, $f'$ in $\Fun(\Delta^o,\C)$ are chain-homotopic in
the sense of Definition~\ref{del.ho.def}, then so are the maps
$E(f)$, $E(f')$ in $\Fun(\Delta^o,\E)$.
\end{lemma}

\proof{} By induction on the length of chain, we may assume that
the maps $f,f':c \to c'$ are connected by a homotopy $\wt{f}:c \to
\Hhom(\Delta_1,c') = \pi_*\pi^*c'$, where $\pi =
\sigma([1])^o:(\Delta/[1])^o \to \Delta^o$ is the projection. But
then we tautologically have $E \circ \pi^* \cong \pi^* \circ E$, and
the adjunction map $a:\pi^*\pi_*\pi^*c' \to \pi^*c'$ gives rise to a
map
$$
E(a):\pi^*E(\pi_*\pi^*c') \cong E(\pi^*\pi_*\pi^*c') \to
E(\pi^*c') \cong \pi^*E(c').
$$
This is in turn adjoint to map $a':E(\pi_*\pi^*c') \to
\pi_*\pi^*E(c')$, and $a' \circ E(\wt{f})$ is an elementary homotopy
connecting $E(f)$ and $E(f')$.
\endproof

If the source category $\C$ is not complete, Lemma~\ref{ch.le} still
holds, with the same proof, if we restrict our attention to a full
subcategory in $\Fun(\Delta^o,\C)$ where $\Hhom(\Delta_1,-)$ and
$\Hhom(\SS_0,-)$ exist. In particular, assume given a small category
$I$ equipped with a covering class $F$; then one can take $\C =
\Cov(i)$, $i \in I$, and consider the subcategory $\HCov(i) \subset
\Fun(\Delta^o,\C)$. As a target category $\E$, let us take a
finitely presentable abelian category. Then for any functor $E:I^o
\to \E$ and object $i \in I$, we can consider the functor $E^o_i =
E^o \circ \sigma(i):\Cov(i) \to \E^o$, and Lemma~\ref{ch.le} applies
to its extension $E^o_i(\Delta):\HCov(i) \to
\Fun(\Delta^o,\E^o)$. Thus if we compose $E^o_i(\Delta)$ with the
Dold-Kan equivalence, we obtain a functor $\HCov(i) \to C_{\geq
  0}(\E^o)$ that sends chain-homotopic map to chain-homotopic maps.
Dually, $E_i:\HCov(i)^o \to \Fun(\Delta,\E) \cong C^{\geq 0}(\E)$
again sends chain-homotopic maps to chain-homotopic maps, so that
for any $n$, the $n$-th homology group $H^n(E_i(i'))$ of the complex
$C^\hdot(E_i(i'))$ defines a functor $\HCov(i)^o \to \E$ that
factors through the filtered category $HC(i)$.

\begin{defn}
The {\em $n$-th Dold-derived functor} $D^n(E):I^o \to \E$ of a
functor $E:I^o \to \E$ is given by
$$
D^n(E)(i) = \colim_{\wt{i} \in HC(i)^o}H^n(C_\idot(E_i(\wt{i}))).
$$
\end{defn}

This is functorial in $i$ since any map $i' \to i$ in $I$ induces a
functor $\HCov(i) \to \HCov(i')$ that sends chain-homotopic maps to
chain-homotopic maps, so we again have a Grothendieck fibration $HC
\to I$ with fibers $HC(i)$, and can apply \eqref{kan.eq}.

\begin{prop}\label{dold.prop}
For any $n \geq 0$, $E \in \Fun(I^o,\E)$, we have a functorial
isomorphism
$$
D^n(E) \cong R^nea(E),
$$
where $a:\Fun(I^o,\E) \to \Shv(I,\E)$ is the associated sheaf
functor, and $R^\hdot e$ are the derived functors of the embedding
functor $e:\Shv(I,\E) \to \Fun(I^o,\E)$.
\end{prop}

\proof{} Take an object $i \in I$, consider the category $\HCov(I)$
with the covering class $F_\Delta$, and choose an $F_\Delta$-hull
$\wt{i} \in \HCov(I)$ that exists by
Lemma~\ref{hull.le}~\thetag{ii}. By definition, $\wt{i}$ is a
proobject in $\HCov(i)$ represented by some filtered diagram
$\gamma:J \to \HCov(i)^o$. Since $\wt{i}$ is $F_\Delta$-liftable,
$i' \setminus J$ is non-empty directed for any $i' \in \HCov(i)$,
and then the same is true for the right comma-fibers of the
composition functor $J \to \HCov(i)^o \to HC(i)^o$. Since $J$ is
filtered, the composition functor is cofinal, and we have
\begin{equation}\label{redu}
D^n(E)(i) \cong \colim_{j \in J}H^n(E_i(\gamma(j))) \cong
H^n(E(\sigma(i)(\wt{i}))),
\end{equation}
where the second isomorphism holds since filtered colimits in $\E$
are exact. Moreover, by Lemma~\ref{ev.le} and adjunction,
$\ev_n(\wt{i})$ is $\sigma(i)^*F$-liftable in $\Cov(i)$, and then
$\sigma(i)(\wt{i})$ is $F$-liftable in $I$ by
Lemma~\ref{hull.le}~\thetag{i}.  By Lemma~\ref{Y.le}, \eqref{redu}
then immediately implies that $D^n(E) \cong D^n(a(E))$, so that
$D^n(-)$ factors through the associated sheaf functor. Moreover, if
$n=0$, Dold-Kand equivalence provides a cartesian square
\begin{equation}\label{H.0}
\begin{CD}
H^0(E_i(i')) @>>> E_i(i'_0)\\
@VVV @VVV\\
E_i(i'_0) @>>> E_i(i'_1)
\end{CD}
\end{equation}
for any hypercovering $i' \in \HCov(i)$, and if $E$ is a sheaf, then
the map $E_i(i'_0 \times_i i'_0) \to E_i(i'_1)$ is injective, so we
can replace $E_i(i'_1)$ in \eqref{H.0} with $E_i(i'_0 \times_i
i'_0)$ and conclude from the sheaf condition that $H^0(E_i(i'))
\cong E(i)$, for any $i' \in \HCov(i')$. Therefore $D^0(E) \cong
e(a(E))$. Finally, by \eqref{redu}, the collection $D^n(-)$ defines a
$\delta$-functor on $\Shv(I,\E)$ in the sense of \cite{toho}, so, if
we say that a sheaf $E \in \Shv(I,\E)$ is {\em exact} if
$H^n(E_i(i')) = 0$ for $n \geq 1$ and any $i \in I$, $i' \in
\HCov(i)$, then it suffices to show that any sheaf $E$ admits an
injective map $E \to E'$ to an exact sheaf $E'$. For any object $i
\in I$, choose an $F$-hull $\wt{i}$ provided by
Lemma~\ref{hull.le}~\thetag{ii}, and consider the product
$$
E' = \prod_{i \in I}\Y_{\wt{i}}(E(\wt{i})),
$$
where $\Y_{\wt{i}}(E(\wt{i}))$, $i \in I$ are as in
Lemma~\ref{Y.le}. Then on one hand, $E'$ is exact by
Lemma~\ref{hur.le}, and on the other hand, the maps $\wt{i} \to i$
induce a map $a:E \to E'$, and since $E$ is a sheaf, all the
maps $E(i) \to E(\wt{i})$ are injective. Thus the map $\ev_i(a)$ is
injective for every $i \in I$, so that $a$ itself is injective.
\endproof

\section{Categories of morphisms.}\label{main.sec}

We can now define and study the main subject of the paper, namely,
the category of morphisms between abelian categories. Assume given
finitely presentable abelian categories $\A$, $\B$. Recall that the
subcategory $\A_c \subset \A$ of compact objects in $\A$ is
preabelian in the sense of Definition~\ref{preab.def}, and the class
$\Epi$ of epimorphisms is a covering class in the opposite category
$\A_c^o$.

\begin{defn}\label{mor.def}
The {\em category of morphisms} $\Mor(\A,\B)$ is the full
subcategory $\Mor(\A,\B) = \Shv(\A_c^o,\B) \subset \Fun(\A_c,\B)
\cong \Fun_c(\A,\B)$ spanned by sheaves with respect to the
$\Epi$-topology on $\A_c^o$.
\end{defn}

By Example~\ref{shv.ab.exa}, the category $\Mor(\A,\B)$ is a
Grothendieck abelian category, and we have a tautological action
functor
\begin{equation}\label{act.mor}
\A \times \Mor(\A,\B) \to \B.
\end{equation}
Conversely, a functor $E:\A \to \B$ comes from a morphism if and
only if
\begin{enumerate}
\item $E$ is continuous, and
\item for any injective map $a:A \to B$ in $\A_c$, the map $E(a):E(A)
  \to E(B)$ is injective in $\B$, and so is the map
\begin{equation}\label{F.sq}
E(B) \oplus_{E(A)} E(B) \to E(B \oplus_A B),
\end{equation}
where we let $B \oplus_A B$ be the cokernel of the map $a \oplus
(-a):A \to B \oplus B$, and similarly for $E(B) \oplus_{E(A)} E(B)$.
\end{enumerate}
Indeed, \eqref{F.sq} is just \eqref{E.sq} for $\langle \A_c^o,\Epi
\rangle$, and we note that epimorphisms in $\A_c^o$ are
monomorphisms in $\A_c$. Since both \thetag{i} and \thetag{ii} are
obviously closed under compositions, we have natural composition
functors
$$
\Mor(\A,\B) \times \Mor(\B,\C) \to \Mor(\A,\C)
$$
for any three finitely presentable abelian categories $\A$, $\B$, $\C$.

\begin{exa}\label{GP0.exa}
The sheaf conditions \thetag{i}, \thetag{ii} immediately imply that
an object $E \in \Fun(\A_c,\B)$ is an additive sheaf if and only if
the its contunuous extension $E:\A \to \B$ is left-exact in the
usual sense (that is, commutes with finite limits). Thus the full
subcategory $\Mor_{add}(\A,\B) \subset \Mor(\A,\B)$ spanned by
additive functors is the category of left-exact continuous functors
$\A \to \B$. In particular, this category of functors is a
Grothendieck abelian category. However, $\Mor(\A,\B)$ is strictly
bigger: not all morphisms are additive. For example, if $\A = \B =
\Ab$, then $\Mor(\A,\B) \subset \Fun_c(\Ab,\Ab)$ contains the
functor sending an abelian group $A$ to the free abelian group
$\Z[A]$ spanned by the underlying set of $A$.
\end{exa}

\begin{exa}\label{GP.exa}
Note that if $\A$ is a small abelian category, then to define the
category $\Mor(\Ind(\A),\Ab) \cong \Shv(\A^o,\Ab)$ with its full
subcategory $\Mor_{add}(\Ind(\A),\Ab) \subset \Mor(\Ind(\A),\Ab)$,
and prove that both are Grothen\-dieck abelian categories, we do not
need to know that $\Ind(\A)$ is abelian. Conversely, the simplest
way to prove that $\Ind(\A)$ is abelian is to show that the Yoneda
embedding \eqref{yo.eq} identifies it with
$\Mor_{add}(\Ind(\A),\Ab)$ -- this is a version of the fundamental
Gabriel-Popescu Theorem, \cite[Ch.\ 5, \S 10]{BD}, and in our
setting, it immediately follows from Example~\ref{ind.exa}.
\end{exa}

On the derived level, we denote by $\DMor(\A,\B)$ the derived
category of the abelian category $\Mor(\A,B)$, and we let
$\DMor^{\leq 0}(\A,\B) \subset \DMor(\A,\B)$ be the full subcategory
spanned by coconnective objects. We have a fully faithful embedding
$R^\hdot e:\DMor^{\geq 0}(\A,\B) \to \D^{\geq 0}(\A_c,\B) \cong
\D^{\geq 0}_c(\A,\B)$ of Proposition~\ref{dold.prop}. We recall that
the target category is obtained by localizing the category of
complexes
$$
C^{\geq 0}(\Fun(\A_c,\B)) \cong
\Fun(\A_c,C^{\geq 0}(\B)) \cong 
\Fun_c(\A,C^{\geq 0}(\B))
$$
with respect to quasiisomorphisms, and one can also use the Dold-Kan
equivalence to replace $C^{\geq 0}(-)$ with $\Fun(\Delta,-)$. This
allows to extend the functor interpretation of $\Mor(\A,\B)$ to
$\DMor^{\geq 0}(\A,\B)$. Namely, for any functor $E:\A_c \to
\Fun(\Delta,\B)$, we have its canonical continuous extension $E:\A
\to \Fun(\Delta,\B)$, and we can further define the {\em Dold
  extension}
\begin{equation}\label{E.dold}
\Dd(E):\Fun(\Delta,\A) \to \Fun(\Delta,\B)
\end{equation}
as the composition
$$
\begin{CD}
\Fun(\Delta,\A) @>{E}>> \Fun(\Delta \times \Delta,\B)
@>{\delta^*}>> \Fun(\Delta,\B),
\end{CD}
$$
where the first functor is $E$ applied termwise, and $\delta:\Delta
\to \Delta \times \Delta$ is the diagonal embedding.

\begin{defn}\label{homo.def}
A continuous functor $C^{\geq 0}(\A) \to C^{\geq 0}(\B)$
is {\em homotopical} if it sends quasiisomorphisms to
quasiisomorphisms.
\end{defn}

\begin{theorem}\label{dold.thm}
A functor $E:\A_c \to C^{\geq 0}(\B) \cong \Fun(\Delta,\B)$
represents an object in $\DMor^{\geq 0}(\A,\B) \subset \D^{\geq
  0}(\A_c,\B)$ if and only if its Dold extension \eqref{E.dold} is
homotopical in the sense of Definition~\ref{homo.def}.
\end{theorem}

\proof{} For the ``only if'' part, say that $E:\A_c \to C^{\geq
  0}(\B)$ is {\em good up to degree $n \geq 0$} if for any
quasiisomorphism $f$ in $C^{\geq 0}(\A)$, $\Dd(E)(f)$ is a
quasiisomorphism in degrees $\leq n$. By Lemma~\ref{DK.le}, a
pointwise quasiisomorphism $E_0 \to E_1$ induces a pointwise
quasiisomorphism $\Dd(E_0) \to \Dd(E_1)$, so that being good only
depends on the object in $\D^{\geq 0}(\A_c,\B)$ represented by
$E$. Moreover, for any complex $M^\hdot \in C^{\geq
  0}(\Mor(\A,\B))$, note that \thetag{i} $R^\hdot e(M^\hdot)$ is
good up to degree $n$ if and only if so is $R^\hdot
e(M^{\leq n+1})$, where $M^{n+1}$ is the stupid truncation $M^0 \to
\dots \to M^{n+1}$, and \thetag{ii} the property of being good up to
degree $n$ is stable under extensions, and the property of being
good up to all degrees is also stable under homological
shifts. Therefore to prove that $R^\hdot e(M^\hdot)$ is good in all
degrees for any complex $M^\hdot$, it suffices to consider complexes
concentrated in degree $0$.

Moreover, say that a sheaf $M \in \Mor(\A,\B)$ is {\em ind-exact} if
$H^n(M(\wt{A})) = 0$ for any $n \geq 1$ and hypercovering $\wt{A}$
of an object $A \in \A^o$ (in particular, this applies to
hypercovering in $\A_c^o \subset \A^o$, so that an ind-exact sheaf
is exact in the sense used in the proof of
Proposition~\ref{dold.prop}). Then since $\A$ is finitely
presentable, it is a Grothendieck abelian category, thus has enough
injectives, and any injective $I \in \A$ is $\Epi$-liftable when
considered as an object of $\A^o \cong \Pro(\A_c^o)$. Then by
Lemma~\ref{Y.le}, $\Y_A(B) \in \Fun(\A_c,\B)$ is a sheaf for any $B
\in \B$, $A \in \A$, and by exactly the same argument as in the
proof of Proposition~\ref{dold.prop}, any $M \in \Mor(\A,\B)$ admits
an embedding into an ind-exact sheaf of the form
$$
M' = \prod_{A \in \A_c}\Y_{\wt{A}}(M(\wt{A})),
$$
where for any $A \in \A_c$, we choose an embedding $A \to \wt{A}$
into some injective $\wt{A} \in \A$. We conclude that every sheaf in
$\Mor(\A,\B)$ admits a resolution by ind-exact sheaves, so it
suffices to prove that $R^\hdot e(M)$ is good for an ind-exact sheaf
$M \in \Mor(\A,\B)$. In this case, the Dold-derived functors
$D^n(M)$ vanish for $n \geq 1$, so we have $R^\hdot e(M) \cong M$ by
Proposition~\ref{dold.prop}, and we have to show that $\Dd(M)(f)$ is
a quasiisomorphism for any quasiisomorphism $f:A^\hdot \to B^\hdot$
in $C^{\geq 0}(\A)$.

In the simple case when $A^\hdot = A \in \A \subset C^{\geq 0}(\A)$
is an object in $\A$ considered as a complex concentrated in
degree $0$, the claim then follows by definition, since $B^\hdot$ is a
hypercovering of the object $A$.

In the general case, note that the abelian category $C^{\geq 0}(\A)$
is also a Grothendieck abelian category, thus has enough injectives,
and moreover, for any injective $I^\hdot \in C^{\geq 0}(\A)$, $I^n
\in \A$ is injective for any $n \geq 0$. Therefore $A^\hdot$ admits
an injective resolution in $C^{\geq 0}(\A)$, and this is a bicomplex
$I^{\hdot,\hdot} \in C^{\geq 0,\geq 0}(\A)$ with injective terms
equipped with a map $A^\hdot \to I^{\hdot,\hdot}$ such that $A^n \to
I^{n,\hdot}$ is a quasiisomorphism for any $n \geq 0$. Then on one
hand, $f_A:A^\hdot \to I^\hdot = \delta^*I^{\hdot,\hdot}$ is a
quasiisomorphism by Lemma~\ref{DK.le}, and $I^\hdot$ is a complex of
injectives, and on the other hand, the simple case of the claim that
we have already proved together with Lemma~\ref{DK.le} show that
$\Dd(M)(f_A)$ is also a quasiisomorphism. Applying the same argument
to $B^\hdot$, we reduce the claim to the case when $f:A^\hdot \to
B^\hdot$ is a map between complexes of injectives. But in this case,
$f$ is invertible up to a chain-homotopy equivalence, and we are
done by Lemma~\ref{ch.le}.

Finally, for the ``if'' part, say that $E \in \D^{\geq 0}(\A_c,\B)$
is good if $\Dd(E)$ is homotopical. Then in particular, $E$ inverts
all maps between $1$-truncated coverings of objects in $\A$ in the
sense of Remark~\ref{trunc.rem}, so that the canonical truncation
$\tau^{\leq 0}E \in \Fun(\A_c,\B) \subset \D^{\geq 0}(\A_c,\B)$ is a
sheaf. Thus if $a(E)=0$, then $\tau^{\leq 0}(E) = 0$, so that $E \in
\D^{\geq 1}(\A_c,\B)$, and then the same argument applies to the
homological shift $E[1]$, so by induction, $E \in \D^{\geq
  n}(\A_c,\B)$ for any $n \geq 0$. We conclude that if $a(E)=0$ for
a good $E$, then $E=0$. If not, let $E'$ be the cone of the
adjunction map $E \to a(R^\hdot e(E))$. Since we have already proved
that $a(R^\hdot e(E))$ is good, so is $E'$, and since $a(E')=0$, we
have $E'=0$ and $E \cong R^\hdot e(a(E))$.
\endproof

By Theorem~\ref{dold.thm}, any continuous functor $E:\A \to
C^\hdot(\B)$ that represents an object in $\DMor^{\geq 0}(\A,\B)$
extends to a homotopical continuous functor $C^{\geq 0}(\A) \to
C^{\geq 0}(\B)$. Let us complement it by showing that any continuous
homotopical functor appears in this way. Namely, recall that
$C^{\geq 0}(\A) \cong \Ind(C^{\geq 0}_b(\A_c))$ is finitely
presentable, and let $\Dh^{\geq 0}(\A,\B) \subset \D_c(C^{\geq
  0}(\A),\B)$ be the full subcategory spanned by homotopical
functors. We then have a natural projection
\begin{equation}\label{tau.eq}
\tau^*:\Dh^{\geq 0}(\A,\B) \to \D^{\geq 0}(\A_c,\B) \cong \D^{\geq
  0}_c(\A,\B),
\end{equation}
where $\tau:\A \to C^{\geq 0}(\A)$ sends $A \in \A$ to itself
considered as a complex concentrated in degree $0$.

\begin{corr}\label{dmor.corr}
The projection \eqref{tau.eq} factors through an equivalence of
categories $\Dh^{\geq 0}(\A,B) \cong \DMor^{\geq 0}(\A,\B) \subset
\D^{\geq 0}(\A_c,\B) \cong \D^{\geq 0}_c(\A,\B)$.
\end{corr}

\proof{} The fact that \eqref{tau.eq} factors through $\DMor^{\geq
  0}(\A,\B)$ immediately follows from
Theorem~\ref{dold.thm}. Moreover, we obviously have $\tau^*(\Dd(E))
\cong E$, so that the correspondence $E \mapsto \Dd(E)$ provides a
one-sided inverse to the functor $\tau^*:\Dh^{\geq 0}(\A,\B) \to
\DMor^{\geq 0}(\A,\B)$. To prove the claim, it then suffices to
construct a functorial isomorphism $E \cong \Dd(\tau^*E)$ for any $E
\in \Dh^{\geq 0}(\A,\B)$. To do this, represent $E$ by a homotopical
continuous functor $E:C^{\geq 0}(\A) \to C^{\geq 0}(\B)$, and note
that $C^{\geq 0}(\A)$ is also an abelian category, so that $E$
admits a Dold extension to a continuous functor $\Dd(E):C^{\geq
  0,\geq 0}(\A) \to C^{\geq 0}(\B)$. We then have $E \cong \Dd(E)
\circ \Ll$ and $\Dd((\tau^*E)) \cong \Dd(E) \circ \Rr$, where $\Ll$
and $\Rr$ are the embeddings of Example~\ref{bimod.exa}. However, we
also have the embedding $\I$ and maps $\Ll,\Rr \to \I$ that give
rise to functorial maps
\begin{equation}\label{LRT.eq}
\begin{CD}
E \cong \Dd(E) \circ \Ll @>>> \Dd(E) \circ \I @<<< \Dd(E) \circ \Rr
\cong \Dd(\tau^*E),
\end{CD}
\end{equation}
so it suffices to prove that both maps in \eqref{LRT.eq} are
quasiisomorphisms. By Example~\ref{bimod.exa}, it then suffices to
prove that $\Dd(E)$ sends both left and right quasiisomorphisms of
Definition~\ref{bimod.def} to quasiisomorphisms in $C^{\geq
  0}(\B)$. For right quasiisomorphisms, this immediately
follows from the fact that $E$ is homotopical, and for left
quasiisomorphisms, this is Theorem~\ref{dold.thm}.
\endproof

\begin{remark}\label{kazh.rem}
As we saw in the proof of Theorem~\ref{dold.thm}, any object $E$ in
$\DMor^{\geq 0}(\A,\B)$ can be represented by a complex $M^\hdot$ of
ind-exact sheaves, and then $E \cong e(M^\hdot):\A_c \to C^{\geq
  0}(\B)$ is also a complex of sheaves, so that the Dold extension
$\D(E)$ is then represented by a continuous homotopical functor
$C^{\geq 0}(\A) \to C^{\geq 0}(\B)$ that sends monomorphisms to
monomorphisms. This makes it into a {\em left-derivable} functor
with respect to the injective model structures on $C^{\geq 0}(\A)$,
$C^{\geq 0}(\B)$ in the sense of \cite[Definition 3.1]{kazh}. By
Corollary~\ref{dmor.corr}, any continuous homotopical functor is
then quasiisomorphic to a left-derivable one. This has its uses; in
particular, left-derivable functors can be used in gluing
constructions such as \cite[Section 3]{kazh}.
\end{remark}

\section{Stability and extensions.}\label{stab.sec}

\subsection{Stability.}

In the assumptions of Theorem~\ref{dold.thm}, any homotopical
functor $E:C^{\geq 0}(\A) \to C^{\geq 0}(\B)$ by definition descends
to a functor
\begin{equation}\label{D.E.C}
\D(E):\D^{\geq 0}(\A) \to \D^{\geq 0}(\B),
\end{equation}
and Theorem~\ref{dold.thm} then shows that \eqref{act.mor} extends
to an action functor
\begin{equation}\label{act.dmor}
\D^{\geq 0}(\A) \times \DMor^{\geq 0}(\A,\B) \to \D^{\geq
  0}(\B).
\end{equation}
However, nothing in Theorem~\ref{dold.thm} insures that the functor
$\D(E)$ is triangulated in any sense of the word. In particular, it
is not necessarily additive, and it does not commute with
holomological shifts, so that there is no easy way to extend
\eqref{act.dmor} to full derived categories. The latter can be
actually made more precise, in the following way.

For any cosimplicial object $c:\Delta \to \C$ in a complete category
$\C$, and any simplicial set $X$, we have the cosimplicial object
$c(X):\Delta \to \C$ with terms $c(X)([n]) = c([n])(X([n]))$;
equivalently, $c(X) \cong \pi_*\pi^*c$, where $\pi:\Delta X \to
\Delta$ is the forgetful functor from the category of elements of
$X:\Delta^o \to \Sets$. If we have a functor $E:\C \to \E$ to
another complete category $\E$, and extend it to functors
$E:\Fun(\Delta,\C) \to \Fun(\Delta,\E)$, $E:\Fun(\Delta X,\C) \to
\Fun(\Delta X,\E)$ by applying it pointwise, then by the same
argument as in Lemma~\ref{ch.le}, the tautological isomorphism $E
\circ \pi^* \cong \pi^* \circ E$ induces by adjunction a map $E
\circ \pi_* \to \pi_* \circ E$, and taken together, the two give a
map
\begin{equation}\label{E.X}
E(c(X)) \to E(c)(X)
\end{equation}
functorial in $E$, $c$ and $X$. Now take $\C = \A$, $\E = \B$, and
let $\SS^1:\Delta^o \to \Sets$ be the simplicial circle obtained by
gluing together both ends of the simplicial interval $\Delta_1$ (this
is different from the standard simplicial circle $\SS_1$ that we
define to be the boundary of the $2$-simplex $\Delta_2$). Then the
projection $\SS^1 \to \Delta_0 = \ppt_{\Delta^o}$ admits a unique
splitting $\Delta_0 \to \SS^1$, so that for any $A:\Delta \to \A$,
we have a splitting
\begin{equation}\label{om.E}
A(\SS^1) \cong A \oplus \Omega(A)
\end{equation}
for some object $\Omega(A):\Delta \to \A$ functorial in $A$.
Up to a quasiisomorphism, $\Omega(A)$ corresponds to the homological
shift $A[-1]$ under the Dold-Kan equivalence \eqref{DK.eq}. The map
\eqref{E.X} then provides a functorial map
\begin{equation}\label{om.st}
E(\Omega(A)) \to \Omega(E(A)),
\end{equation}
and one may ask that this map is a quasiisomorphism for any
$A$. This turns out to be a non-trivial condition on $E$, and it can
be stated in several equivalent ways.

\begin{defn}\label{st.def}
For any two abelian categories $\A$, $\B$, a functor $E:\A \to
C^\hdot(\B)$ is {\em stable} if for any short exact sequence
\eqref{ab1.dia} in $\A$, with the corresponding bicartesian square
$\gamma:[1]^2 \to \A$ of \eqref{ab1.sq}, the induced square
$\gamma^*E \in \D([1]^2,\B)$ is homotopy bicartesian.
\end{defn}

Stability in the sense of Definition~\ref{st.def} is obviously
invariant under cones and quasiisomorphisms, so if $\A$ and $\B$ are
finitely presentable and $E$ is continuous, it only depends on the
object in $\D_c(\A,\B)$ represented by $E$. In particular, it makes
sense to say that an object in $\DMor(\A,\B) \subset \D_c(\A,\B)$ is
stable. Moreover, since $C^{\geq 0}(\A)$ is also abelian and
finitely presentable, stability also makes sense for objects in
$\Dh^{\geq 0}(\A,\B) \subset \D_c(C^{\geq 0}(\A),\B)$.

\begin{remark}
In practice, stability means two things: \thetag{i} the functor
$\D(E)$ of \eqref{D.E} is pointed, so that $\D(E)(p) \circ \D(E)(i)
= \D(E)(p \circ i)$ is $0$, and \thetag{ii} the induced map from a
cone of $\D(E)(i)$ to $\D(E)(C)$ is an isomorphism, so that $\D(E)$
sends short exact sequences to distinguished triangles. However, as
we saw in Example~\ref{D.exa}, this ``induced map'' only becomes
uniquely defined once we lift $\D(E)$ to an object $E \in
\D(\A,\B)$; just knowing $\D(E)$ is not enough.
\end{remark}

\begin{prop}\label{st.prop}
Assume given an object $E \in \D^{\geq 0}_c(\A,\B)$ represented by a
continuous homotopical functor $E:C^{\geq 0}(\A) \to C^{\geq
  0}(\B)$, with the corresponding functor $\D(E)$ of
\eqref{D.E.C}. Then the following conditions are equivalent.
\begin{enumerate}
\item The functor $\D(E)$ is additive.
\item The object $E$ is stable.
\item The object $\tau^*E \in \DMor^{\geq 0}(\A,B)$ is stable.
\item The functor $\D(\tau^*E):\A \to \D(\B)$ is additive.
\item The map \eqref{om.st} is a quasiisomorphism for any $A \in
  C^{\geq 0}(\A)$.
\end{enumerate}
\end{prop}

\proof{} Since stability for split short exact sequences is
equivalent to additivity, \thetag{ii}$\Rightarrow$\thetag{i} and
\thetag{iii}$\Rightarrow$\thetag{iv}. We obviously have
\thetag{iv}$\Rightarrow$\thetag{i} and
\thetag{iii}$\Rightarrow$\thetag{ii} by restriction to $\A \subset
C^{\geq 0}(\A)$. Conversely, since $E$ is quasiisomorphic to the
Dold extension $\Dd(\tau^*E)$ by Corollary~\ref{dmor.corr},
\thetag{iv}$\Rightarrow$\thetag{i} and
\thetag{iii}$\Rightarrow$\thetag{ii}. Moreover, \thetag{iv} actually
implies stability for $E$ and termwise-split short exact sequences
in $C^{\geq 0}(\A) \cong \Fun(\Delta,\A)$, and every
short exact sequence \eqref{ab1.dia} in $\A \subset C^{\geq 0}(\A)$
is quasiisomorphic to the sequence
\begin{equation}\label{ab.C}
\begin{CD}
0 @>>> \Cc(p)[-1] @>>> \Cc(i \oplus i)[-1] @>>> C @>>> 0
\end{CD}
\end{equation}
in $C^{\geq 0}(\A)$, where $i \oplus i:B \oplus_A B \to C$ is the
natural map. Since \eqref{ab.C} is termwise-split,
\thetag{iv}$\Rightarrow$\thetag{iii}, and \thetag{i}, \thetag{ii},
\thetag{iii} and \thetag{iv} are all equivalent. Analogously,
\thetag{iv} applied pointwise shows that the map \eqref{E.X} is a
quasiisomorphism for any finite $X$, thus \thetag{v}, and to finish
the proof, it remains to show that
\thetag{v}$\Rightarrow$\thetag{iv}.

To do this, assume that \thetag{v} holds, and note that since
\eqref{om.st} is a quasiisomorphism for $A=0$, the functor $\D(E)$
is pointed. Thus for any two objects $A_0,A_1 \in \A \subset C^{\geq
  0}(\A)$, with embeddings $i_l:A_l \to A_0 \oplus A_1$, $l=0,1$ and
projections $p_l:A_0 \oplus A_0 \to A_l$, $l=0,1$, the composition
$$
\begin{CD}
E(A_0) \oplus E(A_1) @>{E(i_0) \oplus E(i_1)}>> E(A_0 \oplus A_1)
@>{E(p_0) \oplus E(p_1)}>> E(A_0) \oplus E(A_1)
\end{CD}
$$
is an identity map in $\D(\B)$. Therefore $E \in \D_c(\A \times
\A,\B)$ admits a functorial splitting $E(A_0 \oplus A_1) \cong
E(A_0) \oplus E(A_1) \oplus E'(A_0,A_1)$ for a certain object $E'
\in \D(\A \times \A,\Ind(\B))$, and it suffices to invoke the
following.

\begin{lemma}
Let $\A$, $\B$ be abelian categories, and assume given a functor
$E:\A \times \A \to C^{\geq 0}(\B)$, with the Dold extension
$\Dd(E):C^{\geq 0}(\A \times \A) \to C^{\geq 0}(\B)$ and the
corresponding functor $\D(E):\A \times \A \to \D^{\geq 0}(\B)$ of
\eqref{D.E}. Moreover, assume that for any $A \in \A$, $\D(E)(0
\times A) = \D(E)(A \times 0) = 0$, and for any $A^\hdot \in C^{\geq
  0}(\A \times \A)$, the map
\begin{equation}\label{om.2}
\Dd(E)(\Omega(A^\hdot)) \to \Omega(\Dd(E)(A^\hdot))
\end{equation}
of \eqref{om.st} is a quasiisomorphism. Then $\D(E)=0$.
\end{lemma}

\proof{} Assume that $\D(E) \neq 0$, and let $n$ be the largest
integer such that $\D(E):\A \times \A \to \D^{\geq 0}(\B)$ lands in
$\D^{\geq n}(\B)$. Replacing $E$ with the canonical truncation of
its homological shift $E[n]$, we may assume that $n=0$, and to
derive a contradiction, it suffices to show that $\D(E)$ lands in
$\D^{\geq 1}(\B)$. Since \eqref{om.2} is a quasiisomorphism, it
is further suffices to show that for any $A_0^\hdot,A_1^\hdot \in C^{\geq
  0}(\A)$, $\Dd(E)(\Omega(A_0^\hdot) \times
\Omega(A_1^\hdot)) \in C^{\geq 0}(\B)$ projects into $\D^{\geq
  2}(\B)$. However,
for any $B_0,B_1 \in \Fun(\Delta,\A)$, we have $B_0 \times B_1 \cong
\delta^*(B_0 \boxtimes B_1)$, where the box product $B_0 \boxtimes
B_1 \in \Fun(\Delta \times \Delta,\A \times \A)$ sends $[n] \times
[m] \in \Delta \times \Delta$ to $B_0([n]) \times
B_1([m])$. Therefore $\Dd(E)(B_0 \times B_1)$ can be computed by
taking $B_0 \times B_1 \in \Fun(\Delta,\A) \times \Fun(\Delta,\A)$
and applying the composition
$$
\begin{CD}
\Fun(\Delta,A) \times \Fun(\Delta, A) @>{- \boxtimes -}>>
\Fun(\Delta \times \Delta,\A \times \A)\\
@>{E}>> \Fun(\Delta \times \Delta,C^{\geq 0}(\B))\\
@>{\delta^*}>> \Fun(\Delta,C^{\geq 0}(\B)) \cong \Fun(\Delta \times
\Delta,\B)\\
@>{\delta^*}>> \Fun(\Delta,\B) \cong C^{\geq 0}(\B).
\end{CD}
$$
If we now let $B_0$, $B_1$ be the images of the complexes
$\Omega(A_0^\hdot)$, $\Omega(A_1^\hdot)$ under the Dold-Kan
equivalence, then $B_0([0]) = B_1([0]) = 0$ by the definition of the
functor $\Omega$, and since $E(-,0)$ and $E(0,-)$ are acyclic, the
functor $E(B_0 \boxtimes B_1):\Delta \times \Delta \to C^{\geq
  0}(\B)$ lands in acyclic complexes after restriction to $\Delta
\times [0]$ and $[0] \times \Delta$. In other words, if we apply the
Dold-Kan equivalence \eqref{DDK.eq}, and let $M^{\hdot,\hdot} \in
C^{\geq 0,\geq 0}(C^{\geq 0}(\B))$ be the corresponding bicomplex with
values in $C^{\geq 0}(\B)$, then $M^{n,0}$ and $M^{0,m}$ are acyclic for
any $n,m \geq 0$. This implies that the double totalization
$\Tot(\Tot(M^{\hdot,\hdot})) \in C^{\geq 0}(\B)$ projects into
$\D^{\geq 2}(\B)$, and by virtue of the shuffle quasiisomorphism
\eqref{sh.eq}, the same then holds for $\Dd(E)(\Omega(A_0^\hdot)
\times \Omega(A^\hdot_1)) \cong \delta^*(\delta^*M^{\hdot,\hdot})$.
\endproof

\subsection{Extensions.}

Since stability in the sense of Definition~\ref{st.def} is closed
under taking cones, stable objects form a full triangulated
subcategory $\DMor_{st}(\A,\B) \subset \DMor(\A,\B)$. The standard
$t$-structure on $\DMor(\A,\B)$ induces a $t$-structure on
$\DMor(\A,\B)$, and by Proposition~\ref{st.prop}~\thetag{ii} and
Example~\ref{GP0.exa}, its heart is the category $\Mor_{add}(\A,\B)$
of continuous left-exact functors $\A \to \B$. However, the whole
$\DMor_{st}(\A,\B)$ is {\em not} the derived category of this heart
--- it is bigger. Indeed, Proposition~\ref{st.prop} only imposes
additivity on the derived category level: it is not required, and it
is in general not true, that an object $E$ with additive $\D(E)$ can
be represented by a complex of additive sheaves. It is only the
homology objects of the complex that are required to be additive.

\begin{exa}
Let $\A \cong \B \cong k\amod$ be the category of vector spaces over
a perfect field $k$. Then since $k\amod$ is semisimple, every
left-exact functor $k\amod \to k\amod$ is exact, and if it is also
continuous, then it is completely determined by its value on the
$1$-dimensional vector space $k$, so that $\Mor_{add}(\A,\B) \cong
k\amod$, with $k$ corresponding to the identity functor. However,
$\RHom^\hdot(k,k)$ computed in the category $\DMor_{st}(\A,\B)$ is
Mac Lane Cohomology $HM^\hdot(k)$, and it is highly non-trivial as
soon as $k$ has positive characteristic (in particular, there is a
non-trivial class in $\Ext^2(k,k)$).
\end{exa}

The simplest way to extend \eqref{act.dmor} to full derived
categories -- or at least, to the derived categories $\D^+(-)$ of
complexes bounded from below -- is to use
Proposition~\ref{st.prop}~\thetag{v}: for every stable $E \in
\DMor^{\geq 0}(\A,\B)$, the quasiisomorphism \eqref{om.st} induces a
functorial isomorphism
\begin{equation}\label{sh.st}
\D(E)(A[-1])[1] \cong \D(E)(A),
\end{equation}
and then $\D(E)$ immediately extends to a triangulated functor
\begin{equation}\label{D.E.st}
\D(E):\D^+(\A) = \bigcup_{n \geq 0} \D^{\leq -n}(\A) \to
\D^+(\B)
\end{equation}
by taking the limit with respect to the maps \eqref{sh.st} (the
limit exists since for any individual $A \in \D^+(\Ind(\A))$, the
inverse system stabilizes at a finite step). Since we also have an
obvious identification $E(A)[-1] \cong (E[-1])(A)$, \eqref{act.dmor}
then extends to a functor
\begin{equation}\label{act.dmor.pl}
\D^+(\A) \times \DMor^+_{st}(\A,\B) \to \D^+(\B),
\end{equation}
triangulated separately in each of the two variables. Let us finish
the paper by showing how to lift \eqref{D.E.st} to the chain
level. This is not quite trivial: using \eqref{om.st} directly would
require one to use limits over chain-level liftings of the maps
\eqref{sh.st}, and infinite limits do not behave nicely unless $\B$
satisfies $AB4^*$. Therefore we use an alternative approach based on
``chain-cochain complexes'' as in e.g.\ \cite{trace}.

For any abelian category $\E$, denote by $C^{\geq 0}_{\geq 0}(\E)$
the category of second-quad\-rant bicomplexes in $\E$ (there are
called ``chain-cochain complexes'' in \cite[Section
  3.1]{trace}). Say that a map $f:E^\hdot_\idot \to F^\hdot_\idot$
in $C^{\geq 0}_{\geq 0}(\E)$ is a {\em vertical quasiisomorphism} if
$f:E^\hdot_n \to F^\hdot_n$ is a quasiisomorphism for any $n \geq
0$. We have the sum-totalization functor $\Tot:C^{\geq 0}_{\geq
  0}(\E) \to C^\hdot(\E)$ given by $\Tot(E^\hdot_\idot)^n =
\bigoplus_{i-j=n}E^i_j$. It sends vertical quasiisomorphisms to
quasiisomorphisms and, as in Example~\ref{bimod.exa}, it has a
right-adjoint $\I:C^\hdot(\E) \to C^{\geq 0}_{\leq 0}(\E)$ given by
$\I(E^\hdot)^i_j = E^{i-j} \oplus E^{i-j-1}$. Again as in
Example~\ref{bimod.exa}, we also have full embeddings $\Ll:C_{\geq
  0}(\E) \to C^{\geq 0}_{\geq 0}(\E)$ resp.\ $\Rr:C^{\geq 0}(\E) \to
C^{\geq 0}_{\geq 0}(\E)$ onto the full subcategories of
chain-cochain complexes concentrated in cohomological
resp.\ homological degree $0$, and the isomorphism $\Tot \circ \Rr
\cong \id$ induces a vertical quasiisomorphism $\Rr \to \I \circ
\iota$, where $\iota:C^{\geq 0}(\E) \to C^\hdot(\E)$ is the
tautological embedding. The Dold-Kan equivalence \eqref{DK.eq}
identifies $C^{\geq 0}_{\geq 0}(\E)$ with the category
$\Fun(\Delta^o \times \Delta,\E)$ of simplicial-cosimplicial objects
in $\E$, and a map is a vertical quasiisomorphism if and only if for
any $[n] \in \Delta^o$, it becomes a quasiisomorphism after
restriction to $[n] \times \Delta \subset \Delta^o \times \Delta$.

Now assume given a stable object in $\DMor^{\geq 0}(\A,\B)$
represented by a continuous functor $E:\A \to C^{\geq
  0}(\B)$, consider its Dold extension $\Dd(E)$ of
\eqref{E.dold}, and extend it further to a functor
\begin{equation}\label{E.dd}
\Dd(E):\Fun(\Delta^o \times \Delta,\A) \to \Fun(\Delta^o
\times \Delta,\B)
\end{equation}
by applying it pointwise along $\Delta^o$. We can then consider a
continuous functor
\begin{equation}\label{bid}
\wt{E} = \Tot \circ \Dd(E) \circ \I:C^+(\A) \to C^+(\B),
\end{equation}
where we restrict our attention to complexes bounded from below, and
the map $\Rr \to \I \circ \iota$ provides a map
\begin{equation}\label{R.I}
\iota \circ \Dd(E) \to \wt{E} \circ \iota,
\end{equation}
where $\Dd(E)$ on the left is the Dold extension \eqref{E.dold}.

\begin{prop}
For any continuous functor $E:\A \to C^{\geq 0}(\B)$
representing a stable object in $\DMor^{\geq 0}(\A,\B)$, the functor
\eqref{bid} sends quasiisomorphisms to quasiisomorphisms, and the
map \eqref{R.I} is a quasiisomorphism.
\end{prop}

\proof{} To prove that \eqref{R.I} is a quasiisomorphism, it
suffices to recall that $\Rr \to \I \circ \iota$ is a vertical
quasiisomorphism, and observe that \eqref{E.dd} sends vertical
quasiisomorphisms to vertical quasiisomorphisms by
Theorem~\ref{dold.thm}. For the first claim, we note that every
quasiisomorphism $f:A^\hdot \to B^\hdot$ factors as
\begin{equation}\label{cone}
\begin{CD}
A^\hdot @>{f \oplus i}>> B^\hdot \oplus \Cc(\id_{A^\hdot}) @>{\id \oplus 0}>>
  B^\hdot,
\end{CD}
\end{equation}
where $\Cc(\id_{A^\hdot})$ is the cone of $\id:A^\hdot \to A^\hdot$
and $i:A^\hdot \to \Cc(A^\hdot)$ is the natural embedding; then $f \oplus i$
in \eqref{cone} is an injective quasiisomorphism, and $\id \oplus 0$
admits an injective left-inverse. Therefore it suffices to consider
injective quasiisomorphisms $f:A^\hdot \to B^\hdot$. Moreover, since
$E$ is stable, it suffices to check that $\wt{E}$ sends any acyclic
complex, for example $\Coker f$, to an acyclic complex. We then
observe that for any $n \geq 0$ and acyclic complex $A^\hdot$ in
$C^\hdot(\Ind(\A))$, the complex $\I(A^\hdot)^\hdot_n$ only has
homology in degree $0$, and we have a natural vertical
quasiisomorphism $\Ll(\tau^{\geq 0}A^\hdot) \to \I(A^\hdot)$, where
$\tau^{\geq 0}A^\hdot$ is the canonical truncation. Therefore it
suffices to check that for any acyclic complex $A^\hdot$
concentrated in cohomological degrees $\leq 0$ and bounded from
below, $\wt{E}(A^\hdot)$ is acyclic. However, since $A^\hdot$ is
bounded from below, it has a finite filtration with contractible
associated graded quotients, so since $E$ is stable, it further
suffices to assume that $A^\hdot$ is contractible. But the functor
$\wt{E} \circ \Ll:C^{\leq 0}(\A) \cong C_{\geq 0}(\A)
\to C^\hdot(\Ind(\B))$ is simply the simplicial Dold extension of
the functor $E$, so we are done by Lemma~\ref{ch.le}.
\endproof

\begin{remark}
As in Remark~\ref{kazh.rem}, if a stable object $E \in \DMor^{\geq
  0}(\A,\B)$ is represented by a complex of ind-exact sheaves, then
the resulting functor \eqref{bid} is also left-derivable with
respect to the injective model structures (that is, sends
monomorphisms to monomorphisms).
\end{remark}

{\small\noindent
Affiliations:
\begin{enumerate}
\renewcommand{\labelenumi}{\arabic{enumi}.}
\item Steklov Mathematics Institute (main affiliation).
\item National Research University Higher School of Economics.
\end{enumerate}}

{\small\noindent
{\em E-mail address\/}: {\tt kaledin@mi-ras.ru}
}

\end{document}